\documentclass[letterpaper]{article}
\usepackage[T1]{fontenc}
\usepackage[margin=1in]{geometry}
\usepackage{microtype}
\usepackage{amsmath,amssymb,amsthm,mathtools,booktabs}
\usepackage[authoryear,round]{natbib}
\usepackage[hidelinks]{hyperref}

\theoremstyle{plain}
\newtheorem{theorem}{Theorem}[section]
\newtheorem{proposition}[theorem]{Proposition}
\newtheorem{lemma}[theorem]{Lemma}

\theoremstyle{definition}

\providecommand{\MR}[1]{}
\newcommand{\ARXIV}[1]{\href{https://arxiv.org/abs/#1}{arXiv:#1}}
\newcommand{\E}{\mathbf E}
\newcommand{\PP}{\mathbf P}
\newcommand{\N}{\mathbb N}
\newcommand{\bc}{\overline c}
\newcommand{\one}{\boldsymbol 1}
\newcommand{\diag}{\operatorname{diag}}

\newcommand{\Law}{\mathcal L}
\newcommand{\TV}{\mathrm{TV}}
\newcommand{\dd}{\mathrm d}

\DeclareRobustCommand{\routingd}{d}
\hypersetup{
  pdftitle={Critical heavy traffic and diffusion limits for shortest-of-d dynamic routing},
  pdfauthor={Elmira Yu. Kalimulina}
}
\title{Critical heavy traffic and diffusion limits\\ for shortest-of-\texorpdfstring{$\routingd$}{d} dynamic routing}
\author{Elmira Yu.\ Kalimulina\\
  Lomonosov Moscow State University, Moscow, Russia\\
  \href{mailto:elmira.yu.k@gmail.com}{\texttt{elmira.yu.k@gmail.com}}}
\date{}

\begin{document}

\maketitle
\begin{abstract}
We study shortest-of-$d$ routing, with fixed $d\ge2$, in $N$ parallel
queues with independent finite-state Markov service environments switching
at rate $\kappa_N$. At the critical load $\rho_N=1-1/N$, the empirical measure of the centered queue lengths, under stationarity,
converges in distribution to a random probability measure determined by its mean. We identify the
mean's density and the conditional independence of fixed finite collections
of centered queues in the limit. On time scales of order $N$, the centered empirical mean
converges to a diffusion with explicit drift and constant diffusion
coefficient. These limits hold when $\kappa_N/N\to\infty$, along
subsequences with convergent lattice phase. The diffusion limit extends to
centered initial empirical measures tight in the first Wasserstein topology, with means
converging in law. A short initial layer prepares the limiting family
while preserving the limiting mean. Empirical approximation holds in time-integrated
Wasserstein distance.
\end{abstract}

\noindent\textbf{MSC 2020.} Primary 60K25, 60F17. Secondary 60J27, 60K37.

\noindent\textbf{Keywords.} Heavy traffic; Dynamic routing; Conditional independence; Stationary limits; Diffusion approximation; Markov environment.

\section{Introduction}
This work concerns dynamic routing in parallel queueing systems and is motivated
by the shortest-of-two model of \citet*{VDK}. Each arriving customer samples
two servers uniformly and joins the shorter queue. Its shortest-of-$d$
generalization, with fixed $d\ge2$, is also known as the \emph{supermarket
model} \citep{Mitzenmacher2001}. For Poisson arrivals and exponential service
at a common constant rate, the stationary mean-field queue-length distribution $\zeta_\rho$ has tail probabilities
\begin{equation}\label{eq:classical-tail}
 z_k:=\zeta_\rho([k,\infty))=\rho^{(d^k-1)/(d-1)},\qquad k\in\N_0,
\end{equation}
at each fixed load $0<\rho<1$. We ask what replaces this deterministic
description when the load approaches one as the number of queues grows.

Our main result describes the random limiting empirical measure through
a single scalar coordinate, its mean, at the critical gap
$1-\rho_N=N^{-1}$, where total spare capacity is of order one. At stationarity, the empirical measure of the centered
queue lengths $Q_i-m_N$, where $m_N=\lfloor\log_dN\rfloor$, converges in
distribution to a random probability measure $\Pi_{X_\theta}$. The family
$\Pi_x$ has double-exponential right tails and is parameterized by its mean $x$.
Theorem~\ref{th:critical-random-profile} identifies the density of
$X_\theta$; in the limiting joint law, any fixed number of centered
queues are independent conditional on $X_\theta$, with common law
$\Pi_{X_\theta}$. Limits are taken along subsequences on which the
fractional part of $\log_dN$ converges to $\theta$.

The centered empirical mean is the slow coordinate of this family.
Theorem~\ref{th:initial-layer-diffusion} gives, for every finite $T$,
\[
 \left(\frac1N\sum_i Q_i(Nt/\bc)-m_N\right)_{0\le t\le T}
 \ \Longrightarrow\ (X_t)_{0\le t\le T}
 \quad\text{in }D([0,T],\mathbb R)\text{ with }J_1,
\]
where $\bc$ is the mean service capacity. The limit is the nonlinear
Langevin diffusion $\dd X_t=-V_\theta'(X_t)\,\dd t+\sqrt2\,\dd B_t$,
with the explicit confining potential $V_\theta$ defined in
Section~\ref{sec:model} and the non-Gaussian invariant density
$p_\theta\propto e^{-V_\theta}$ identified in the stationary limit.
The theorem allows centered initial empirical measures tight in the
first Wasserstein topology and centered means converging in distribution.
An initial layer brings the empirical measure close to this family
without changing the limiting mean. Thus no initial approximation by
$\Pi_x$ is required: the mean-process limit holds from time zero, and
empirical approximation is integrated over time.

Each arrival or departure changes the empirical mean by $1/N$, and the
total arrival and departure rate is of order $N$. Accelerating time by $N$ therefore
retains stochastic fluctuations. The key difficulty is the drift. For constant service, it is the number
of idle servers minus one on the accelerated clock. Empty-queue
probabilities must therefore be resolved at the $1/N$ scale, beyond
ordinary convergence of empirical measures. Weighted and time-integrated
balance identities yield the effective drift $g_\theta-1$, rather than
a pointwise approximation of the idle count.

Independently driven finite-state service environments, including phases of zero
capacity, are allowed. Their switching rates are multiplied by $\kappa_N$,
and the critical limits hold under $\kappa_N/N\to\infty$. The number of
choices and the unit-speed environment remain fixed. For constant service,
the queue-only conclusions require no switching parameter. Thus fast service fluctuations preserve the critical limit. The switching
condition is sufficient; its sharpness is not asserted.

\subsection{Related work}
At fixed subcritical load, \citet{LM2007} establish quantitative propagation
of chaos, \citet{Graham2005} proves a stationary functional central limit theorem
for the empirical tail vector, with an infinite-dimensional
Ornstein--Uhlenbeck limit,
and \citet{BF2019} prove finite-horizon central limit theorems for occupancy
measures. Independent service interruptions are studied by \citet{VS2005},
and fluctuating capacities by \citet{JDD2026}. Our critical regime retains
order-one randomness in the logarithmically centered empirical law.

For unbounded power-of-two routing, \citet[Theorem~4.1]{Ying2017} gives
quantitative approximation by a deterministic stationary law when
$\rho_N=1-\gamma N^{-\alpha}>3/4$, $0<\gamma<1$ and $0\le\alpha<1/5$.
Our results identify the random stationary limit and its dynamics at the
critical gap $N^{-1}$.

The scale $N(1-\rho_N)=O(1)$ is suggested by the fixed-load analysis and
numerical conjecture of \citet{Gast2017}; \citet[Theorem~5]{EG2016} provide
a stationary lower bound. At order-one spare capacity, \citet{GW2019}
obtain a full-JSQ diffusion for an uncentered mean of order one.
\citet{BFL2018} study equilibrium concentration when the arrival rate
tends to one and the number of choices grows. \citet{BBD2022} study
near-equilibrium occupancy diffusions with
$d_N\to\infty$, including a regime with exponential drift. For fixed $d$,
\citet{HLM2022} obtain an exponential stationary limit at
$1-\rho_N=N^{-\alpha}$, $\alpha>3$, after multiplying the mean by
$N^{1-\alpha}$. The fixed-$d$ critical regime studied here instead has
logarithmic centering and a nondegenerate random empirical law.

Recent related limits concern a small fraction of choice-routed arrivals
against a random-routing background \citep{AW}, modulated full JSQ at
fixed system size \citep{HLG2026}, and growing choices in the
sub-Halfin--Whitt regime \citep{VCM2025}. \citet{BBL2026} obtain a reflected
infinite Atlas limit for marginal JSQ with a random-routing background;
their stationary convergence result concerns a modified system with pauses.
These works address different routing mechanisms or asymptotic scales.

The proof uses entropy estimates, established averaging methods
\citep{BMP,LN}, lattice asymptotics \citep{Janson2006}, nonlinear stability
ideas \citep{OK2000}, and nonexpansive coupling
\citep[Lemma~2.3]{LM2006}. Section~\ref{subsec:proof-strategy} explains how
these tools enter the critical argument and the initial-layer construction.

Section~\ref{sec:model} states the main results and proof strategy;
Sections~\ref{sec:critical-random-profile}--\ref{sec:initial-layer} prove them.
Technical estimates are proved in the supplement.

\section{Model, main results and proof strategy}\label{sec:model}
\subsection{Model and notation}
Fix an integer $d\ge2$. For $N\ge1$, $0<\rho<1$, and $\kappa>0$, consider
$N$ single-server queues. Let $Q_i(s)$ denote the number of customers at
server $i$, including the customer in service, and let $J_i(s)$ denote its
environment phase; $s$ denotes physical time. A Poisson stream of rate $N\lambda$ feeds the system.
Each arrival samples $d$ indices independently and uniformly, with
replacement, and joins a shortest sampled queue, breaking ties uniformly.
Customers are never discarded.

The phase space is $E=\{1,\ldots,m\}$. Let $A=(a_{ef})$ be an irreducible
unit-speed generator, where $a_{ef}$ is the transition rate from $e$ to $f$,
and let $\pi$ be its invariant probability vector. We regard $\pi$ as a
column vector, so $A^\top\pi=0$ and $\one^\top\pi=1$. Each $J_i$ evolves
autonomously with generator $\kappa A$, including when its queue is empty.
The phase-transition mechanisms are independent across servers. The
randomness driving arrivals, routing, service, and phase transitions is
mutually independent and independent of the initial state; the initial
queue lengths and phases may have an arbitrary joint distribution.

For a capacity function $c:E\to[0,\infty)$, represented as a column vector, put
\[
 \begin{gathered}
 c_* =\max_e c_e,\qquad \bc=\pi^\top c>0,\qquad \lambda=\bc\rho.
 \end{gathered}
\]
Thus $\rho$ is the offered load relative to the average service capacity
$\bc$. The parameters $A,c,d$ are fixed, and constant capacity is permitted.
When a queue is nonempty and its phase is $e$, its length decreases by one
at rate $c_e$. These rates correspond to FCFS service with independent
unit-mean exponential work requirements processed at the current capacity.
Conditional on the initial queue lengths and phases, the remaining work
requirements of all initially present customers are independent unit-mean
exponential random variables. Zero capacity pauses service. It is the required
work, rather than
the calendar service time in the changing environment, that is exponentially
distributed.

For a queue-length configuration $x\in\N_0^N$, write
$U_k(x)=N^{-1}\sum_i\one_{\{x_i\ge k\}}$ and $U_0=1$.
These are the empirical tail probabilities of the finite system. The birth
intensity at server $i$ is
\begin{equation}\label{eq:finite-birth-model}
 b_i(x)=\lambda S_d(U_{x_i},U_{x_i+1}),\qquad
 S_d(u,v)=\sum_{r=0}^{d-1}u^{d-1-r}v^r.
\end{equation}
Indeed, the sampled minimum is $k$ with probability $U_k^d-U_{k+1}^d$.
By symmetry, arrivals with this minimum are shared equally among the
$N(U_k-U_{k+1})$ queues of length $k$, yielding~\eqref{eq:finite-birth-model}.
Hence $b_i\le d\lambda$ and $\sum_i b_i=N\lambda$. Together with service
and phase transitions, these rates specify the Markov chain. Write
$\mu_{N,\rho,\kappa}$ for its invariant law (existence proved below), and
$\mu_{N,\rho,\kappa}^{[\ell]}$ for the marginal on $\ell$ tagged queue--phase pairs. We use $d_{\TV}(P,Q)=\sup_B|P(B)-Q(B)|$.

The finite system has a unique exchangeable invariant probability law;
its phase marginal is $\pi^{\otimes N}$, and its queue lengths have finite
first moments (Lemma~\ref{lem:finite-recurrence}).

\subsection{The critical family and stationary limit}
In the critical regime $\rho_N=1-N^{-1}$, use throughout
\begin{equation}\label{eq:critical-centering}
 m_N=\lfloor\log_dN\rfloor,\quad
 \theta_N=\log_dN-m_N,\quad
 M_N=N^{-1}\sum_iQ_i,\quad X_N=M_N-m_N.
\end{equation}
After centering, the lower boundary recedes to minus infinity. The critical
tail recursion $u_{j+1}=u_j^d$ admits the following family on $\mathbb Z$;
its mean selects the parameter. For $\omega>0$, let $\xi_\omega$ be the law with tails
\begin{equation}\label{eq:critical-profile-family}
 \xi_\omega([j,\infty))=e^{-\omega d^j},\qquad j\in\mathbb Z.
\end{equation}
Its mean
\begin{equation}\label{eq:critical-profile-mean}
 r(\omega)=\sum_{j\ge1}e^{-\omega d^j}
             -\sum_{j\le0}(1-e^{-\omega d^j})
\end{equation}
is a smooth decreasing bijection from $(0,\infty)$ onto $\mathbb R$.
Write $\omega(x)=r^{-1}(x)$ and $\Pi_x=\xi_{\omega(x)}$, the member of this
family with mean $x$. The inverse satisfies
\begin{equation}\label{eq:critical-inverse-mean}
 \omega(x+1)=\omega(x)/d,\qquad
 C_1d^{-x}\le\omega(x)\le C_2d^{-x}.
\end{equation}
The mean is an intrinsic coordinate of this family:
\[
 W_1(\Pi_x,\Pi_y)=|x-y|,\qquad x,y\in\mathbb R,
\]
where $W_1$ is the first Wasserstein distance for cost $|j-j'|$ on
$\mathbb Z$ (Lemma~\ref{lem:nonstationary-W1-input}). Thus motion of the
limiting empirical measure along this family is measured exactly by
motion of its mean.
Constants may depend on the fixed model parameters. For $\theta\in[0,1]$, put
\begin{equation}\label{eq:critical-mixing-density}
 \begin{gathered}
 g_\theta(x)=(d-1)d^\theta\omega(x),\qquad
 V_\theta(x)=x+\int_x^\infty g_\theta(v)\,\dd v,\\
 p_\theta(x)=\mathcal Z_\theta^{-1}e^{-V_\theta(x)},\qquad x\in\mathbb R.
 \end{gathered}
\end{equation}
The integral is finite and $0<\mathcal Z_\theta<\infty$. This density has an
exponential right tail and a faster than exponential left tail; $\theta=1$
denotes a one-sided phase limit.
The function $g_\theta$ represents the effective number of idle servers
after averaging. Ordinary $o(1)$ total-variation error cannot resolve
empty-queue probabilities of order $1/N$. Identifying $g_\theta$ requires
weighted and time-integrated balance; it does not assert pointwise
approximation of the idle count.

The next theorem identifies the random stationary empirical measure through its
actual empirical mean, and hence the residual dependence between tagged queues.
\begin{theorem}[Critical stationary mixture]
\label{th:critical-random-profile}
Suppose $N\to\infty$ and
\[
 \rho_N=1-N^{-1},\qquad \kappa_N/N\longrightarrow\infty,
 \qquad \theta_N\longrightarrow\theta\in[0,1].
\]
Under the finite-system invariant law,
\begin{equation}\label{eq:critical-mean-limit}
 X_N\ \Longrightarrow\ X_\theta,
 \qquad \Law(X_\theta)(\dd x)=p_\theta(x)\,\dd x.
\end{equation}
Let
$\widetilde\mu_N=N^{-1}\sum_i\delta_{(Q_i-m_N,J_i)}$.
Then the following approximation uses the actual empirical mean:
\begin{equation}\label{eq:critical-random-empirical-TV}
 \E d_{\TV}\bigl(\widetilde\mu_N,
                        \Pi_{X_N}\otimes\pi\bigr)\longrightarrow0.
\end{equation}
Consequently $\widetilde\mu_N$ converges in distribution, in the
total-variation metric on probability measures on $\mathbb Z\times E$,
to the random law $\Pi_{X_\theta}\otimes\pi$.
For every fixed $\ell$, the law of the centered designated
queue--phase pairs converges in total variation to
\begin{equation}\label{eq:critical-tagged-mixture}
 \int_{\mathbb R}
       (\Pi_x\otimes\pi)^{\otimes\ell}
                      p_\theta(x)\,\dd x.
\end{equation}
For constant service capacity, the queue-only assertions hold without
any switching parameter or condition on it.
\end{theorem}

\subsection{Diffusion and general initial states}
We distinguish physical time $s$ from accelerated time $t$ by
$s=Nt/\bc$. All accelerated empirical processes below use this clock:
\begin{equation}\label{eq:main-empirical-clocks}
 \begin{aligned}
 P_N(s)&=\frac1N\sum_i\delta_{Q_i(s)-m_N},
 &\mu_N^q(t)&=P_N(Nt/\bc),\\
 \mathcal X_N(t)&=\frac1N\sum_iQ_i(Nt/\bc)-m_N,
 &\widetilde\mu_N(t)&=\frac1N\sum_i
       \delta_{(Q_i(Nt/\bc)-m_N,J_i(Nt/\bc))}.
 \end{aligned}
\end{equation}
At stationarity, $X_N$ and $\widetilde\mu_N$ without a time argument denote
one-time versions. Let $\mathcal P_1(\mathbb Z)$ be the probability measures
with finite first moment, equipped with $W_1$. The candidate limiting mean is the diffusion
\begin{equation}\label{eq:critical-diffusion-sde}
 \dd X(t)=(g_\theta(X(t))-1)\,\dd t+\sqrt2\,\dd B(t).
\end{equation}
Here $B$ is standard Brownian motion, independent of the initial value
specified below. The solution is unique and nonexplosive.
Since $V_\theta'=1-g_\theta$, this is the reversible Langevin diffusion
in potential $V_\theta$, with invariant density $p_\theta$. The potential
is strictly convex: $r'<0$ implies $V_\theta''=-g_\theta'>0$.
For measure-valued paths we use the metric
$\int_0^T d_{\TV}(\nu_t,\nu'_t)\,\dd t$, identifying paths that agree
almost everywhere in time.

\begin{theorem}[Critical diffusion and initial-layer relaxation]
\label{th:initial-layer-diffusion}\label{cor:initial-layer-profile}
Suppose $\rho_N=1-N^{-1}$, $\kappa_N/N\to\infty$ and
$\theta_N\to\theta\in[0,1]$. Assume
\begin{equation}\label{eq:initial-layer-initial-mean}
 Z_N=\frac1N\sum_iQ_i(0)-m_N\ \Longrightarrow\ \nu
\end{equation}
for a probability law $\nu$ on $\mathbb R$, and that the random initial
empirical measures $P_N(0)$ are tight in $\mathcal P_1(\mathbb Z)$.
Equivalently, the latter condition can be written as
\begin{equation}\label{eq:initial-layer-tail-condition}
 \lim_{R\to\infty}\limsup_{N\to\infty}
 \PP\left(\frac1N\sum_i |Q_i(0)-m_N|
       \one_{\{|Q_i(0)-m_N|>R\}}>\eta\right)=0,
 \qquad \eta>0.
\end{equation}
No exchangeability or initial phase restriction is imposed.
The following conclusions hold.
\begin{enumerate}
\item[(i)] For every deterministic physical-time sequence satisfying
\begin{equation}\label{eq:initial-layer-time-choice}
 b_N\longrightarrow\infty,\qquad b_N=o(\log N),
\end{equation}
the initial layer preserves the mean and prepares the empirical measure:
\begin{align}
 &\sup_{0\le s\le b_N}
 \left|\frac1N\sum_iQ_i(s)-\frac1N\sum_iQ_i(0)\right|
                 \longrightarrow0\quad\hbox{in probability},
                                      \label{eq:initial-layer-main-mean}\\
 &W_1\left(P_N(b_N),
       \Pi_{Z_N}\right)
                 \longrightarrow0\quad\hbox{in probability}.
                                      \label{eq:initial-layer-preparation}
\end{align}
\item[(ii)] For each such $b_N$, both the accelerated mean $\mathcal X_N(t)$
started at time zero and the
process observed after the initial layer,
\begin{equation}\label{eq:initial-layer-shifted-process}
 \mathcal X_N^{\rm lay}(t)
     =\frac1N\sum_iQ_i(b_N+Nt/\bc)-m_N,
\end{equation}
converge in $D([0,T],\mathbb R)$ with the $J_1$ topology, for every
finite $T$, to the diffusion~\eqref{eq:critical-diffusion-sde} with $X(0)\sim\nu$,
where $B$ is independent of $X(0)$. No moment assumption on $\nu$ is
needed.
\item[(iii)] For every finite $T$, empirical approximation holds from
accelerated time zero in the following integrated forms:
\begin{align}
 &\int_0^T W_1\left(\mu_N^q(t),\Pi_{\mathcal X_N(t)}\right)\,\dd t
                    \longrightarrow0\quad\hbox{in probability},
                                      \label{eq:initial-layer-time-W1}\\
 &\E\int_0^T d_{\TV}\left(\widetilde\mu_N(t),
                              \Pi_{\mathcal X_N(t)}\otimes\pi\right)\,\dd t
                    \longrightarrow0.\label{eq:initial-layer-time-TV}
\end{align}
The pair $(\mathcal X_N,\widetilde\mu_N)$ converges jointly to
$(X,\Pi_X\otimes\pi)$ in the product of $J_1$ and the integrated
total-variation metric defined above.
\end{enumerate}
For constant service, the queue-only versions need no switching condition.
\end{theorem}

One explicit choice is $b_N=\bc^{-1}\log\log(N+e)$.
Since $\bc b_N/N\to0$, the mean diffusion starts at accelerated time zero;
the empirical measure's initial mismatch is compatible with the integrated topology.
The initial condition controls first-moment tails in probability, without
requiring expectations over the initial randomness. It follows from
tightness of the centered means and, for some $\eta>0$,
\[
 \frac1N\sum_i
 \exp\left(\eta\left|Q_i(0)-\frac1N\sum_lQ_l(0)\right|\right)
       =O_{\PP}(1).
\]
Together, these bounds control first-moment tails on events of
arbitrarily high probability.
Equal queues of length $m_N$ satisfy the condition, although
$\delta_0\ne\Pi_0$. The tail condition prevents a vanishing fraction of queues from carrying
a nonvanishing part of the centered first moment.

The scalar diffusion therefore tracks the limiting empirical measure in
its exact Wasserstein coordinate. At stationarity, mixing the conditional
double-exponential tails gives an exponential unconditional queue-length
tail, derived at the
end of Section~\ref{sec:critical-random-profile}. An equivalent logistic-type
coordinate is recorded in Supplement, Section~S6.

\subsection{Proof strategy}\label{subsec:proof-strategy}
The proof separates the fast adjustment of the empirical measure from
the slow motion of its mean. All estimates used in the stationary argument
hold directly at the critical load; the two main theorems are proved in
the following order.
\begin{enumerate}
\item \emph{Determine the distribution from its mean.}
The finite-system inputs in Proposition~\ref{prop:critical-inputs} control
the weighted square error in the recursion $U_{k+1}\simeq\rho_N U_k^d$.
At an intermediate level $k_N=\lfloor\frac14\log_dN\rfloor$, put
$D_l=N(1-U_l)$ and $Y_N=D_{k_N}/d^{k_N}$. Iterating the recursion determines
the centered empirical tails from $Y_N$. A uniform spread bound then
identifies $Y_N$ with $d^{\theta_N}\omega(X_N)$, including the weighted
error estimates needed for stationary balance.
\item \emph{Identify the stationary law.}
Corrected stationary generator identities give the drift $g_\theta-1$
and identify the mean's density $p_\theta$. Empirical approximation and
exchangeability give the conditional independence in
Theorem~\ref{th:critical-random-profile}. These two steps form
Section~\ref{sec:critical-random-profile}.
\item \emph{Prove the diffusion limit.}
Section~\ref{sec:diffusion-and-starts} first proves the process limit
in stationarity. Uniform integrability and integrated balance replace
accumulated idleness by the amplitude, and the corrected mean martingale
has quadratic variation converging to $2t$. Coupling with reference systems
initialized from invariant laws conditioned on short intervals of the
centered mean transfers this limit to prepared initial states,
with arbitrary initial phases.
\item \emph{Prepare general states.}
The averaged whole-line evolution conserves its mean and relaxes in $W_1$
to $\Pi_x$. Section~\ref{sec:initial-layer} approximates the finite system
on growing physical-time intervals and transfers this relaxation to
$b_N=o(\log N)$. The mean changes negligibly during this layer. Restarting
from the resulting prepared state proves
Theorem~\ref{th:initial-layer-diffusion}.
\end{enumerate}

Table~\ref{tab:main-notation} summarizes the notation.
\begin{table}[t]
\caption{Main notation. Physical time $s$ and accelerated time $t$ satisfy $s=Nt/\bc$.}
\label{tab:main-notation}
\centering\begin{tabular}{@{}p{0.26\linewidth}p{0.69\linewidth}@{}}
\toprule
Notation & Meaning \\
\midrule
$c_e,\bc,\rho,\kappa$ & Service capacity, its stationary mean, load and switching rate.\\
$U_k,D_k$ & Empirical tail and number of queues below level $k$.\\
$m_N,\theta_N$ & Integer and fractional parts of $\log_dN$.\\
$\xi_\omega,r(\omega),\omega(x),\Pi_x$ & Queue-length family, its mean, inverse mean map, and member with mean $x$.\\
$X_N,\mathcal X_N(t)$ & Stationary centered mean and accelerated mean process.\\
$P_N(s),\mu_N^q(t),\widetilde\mu_N(t)$ & Centered empirical queue measures in physical/accelerated time, and joint queue--phase empirical measure.\\
$Y_N,\mathcal R_N,H_N$ & Intermediate amplitude, weighted recursion error, and exponential spread.\\
$g_\theta-1,V_\theta,p_\theta$ & Diffusion drift, potential and invariant density.\\
$S_s,v_j(\mu),V_j^N(s)$ & Whole-line flow, its tails, and centered finite-system tails.\\
\bottomrule\end{tabular}
\end{table}

\section{Critical stationary limit}\label{sec:critical-random-profile}
We prove Theorem~\ref{th:critical-random-profile} by linking the shape of
the centered empirical queue-length measure to the boundary behaviour of
the finite system. An intermediate-level amplitude determines the central
tail probabilities, while empty queues determine the drift of the empirical
mean. The main step connects the averaged contribution of empty queues to
this amplitude and then to the mean. The resulting stationary balance
identity identifies the mixing density.

\subsection{Finite-system estimates}
We first collect the estimates that justify stationary tests, control the
error in the tail recursion, and relate central tails to the empirical mean.
Throughout this subsection expectations are under the finite-system
invariant law, and constants depend only on the fixed parameters unless
stated otherwise. Full proofs of Lemma~\ref{lem:finite-recurrence},
Proposition~\ref{prop:critical-inputs}, and
Lemma~\ref{lem:boundary-spread} are in Supplement, Sections~S1, S2, and~S3,
respectively; we give the main mechanisms here.

To average the service environment, put $\boldsymbol g=c-\bc\one$.
Since $\pi^\top\boldsymbol g=0$, there is a unique centered solution $h$
of the Poisson equation
\begin{equation}\label{eq:poisson}
 -Ah=\boldsymbol g,\qquad \pi^\top h=0.
\end{equation}
For the environment with generator $\kappa A$, the corresponding corrector
is $h/\kappa$. The following Lyapunov estimate justifies the stationary
identities used below.

\begin{lemma}[Finite-system recurrence]\label{lem:finite-recurrence}
For every $N\ge1$, $0<\rho<1$ and $\kappa>0$, the finite system has a
unique invariant probability law, which is exchangeable. For every fixed
$\kappa_*>0$ there is $B<\infty$, depending only on $A,c,d,\kappa_*$, such that
\begin{equation}\label{eq:finite-first-moment}
 \E_{\mu_{N,\rho,\kappa}}Q_i\le B/(1-\rho),\qquad \kappa\ge\kappa_*.
\end{equation}
Its phase marginal is $\pi^{\otimes N}$.
\end{lemma}
\begin{proof}[Proof sketch]\label{prf:lem:finite-recurrence}
Fix $\kappa_*>0$ and $\kappa\ge\kappa_*$. Put
$r_e=2(\|h\|_\infty-h_e)/\kappa\ge0$ and
$V(q,e)=\sum_i(q_i^2+r_{e_i}q_i)$. The routing identity
$\sum_iq_i b_i=N\lambda\sum_{k\ge1}U_k^d\le\lambda\sum_iq_i$
and~\eqref{eq:poisson} give
\[
 \mathcal LV\le-2\bc(1-\rho)\sum_iq_i+CN,
\]
where $C$ depends only on $A,c,d,\kappa_*$. Capping queues and rejecting
arrivals to full queues preserves this bound, since the rejected
increments of $V$ are nonnegative. The capped invariant laws are tight
for fixed $N$; their balance equations pass to an invariant law of the
uncapped chain, and Fatou's lemma preserves the first-moment bound.
Irreducibility gives uniqueness, permutation symmetry gives
exchangeability, and autonomy of the phases gives their marginal
$\pi^{\otimes N}$. For any fixed $\kappa>0$, take $\kappa_*=\kappa$.
\end{proof}

Let $\mathcal L^0$ be the generator with every service capacity replaced by
$\bc$, and let $\mathcal L_{\rm q}$ be the queue-jump part of the actual
generator. For a queue-only test function $F$, put
\[
 S_i^F(x)=\one_{\{x_i>0\}}[F(x-e_i)-F(x)],\qquad
 B_F(x,e)=\sum_i h_{e_i}S_i^F(x).
\]
The phase generator applied to $B_F$ cancels the difference between the two
service drifts, giving the exact identity
\begin{equation}\label{eq:critical-corrected-test}
 \mathcal L(F+B_F/\kappa)
       =\mathcal L^0F+\kappa^{-1}\mathcal L_{\rm q}B_F.
\end{equation}
The stationary balance tests below are bounded or integrable with bounded
jump increments, justifying compensation by the bounded total jump rate
for fixed $N,\kappa$. Exponential Lyapunov tests are handled by truncation.
\begin{proposition}[Critical recursion and deficit bounds]
\label{prop:critical-inputs}
Suppose $N\ge4$, $\rho_N=1-N^{-1}$ and $\kappa_N\ge N$. Put
\begin{equation}\label{eq:critical-residual}
 \mathcal R_N=\sum_{k\ge0}d^{-k-1}
                 (\rho_N U_k^d-U_{k+1})^2,\qquad D_l=N(1-U_l).
\end{equation}
There is a constant $C$, depending only on $A,c,d$, such that
\begin{equation}\label{eq:critical-input-bounds}
 \E D_l\le C\min\{N,d^l\},\qquad
 \E D_l^2\le Cd^{2l}(\log N)^2,\qquad
 \E\mathcal R_N\le C(\log N)^2/N^2,
 \quad l\ge0.
\end{equation}
For constant service the same bounds hold without a switching condition.
\end{proposition}

For $\eta=N^{-3}$ and $z_k=\rho_N^{(d^k-1)/(d-1)}$, use the smoothed entropy
\[
 V_\eta(U)=\sum_{k\ge1}d^{-k}\int_0^{U_k}
                   \log\frac{v+\eta}{z_k}\,\dd v.
\]
Its nonnegative dissipation $\mathcal D_\eta$, defined in Supplement,
Section~S2, controls the residual through
$\mathcal R_N\le C\mathcal D_\eta+C\eta^2$.
The condition $\kappa_N\ge N$ gives
$b=N^{-1}+\kappa_N^{-1}\le2/N$. Corrected stationary balance then yields,
with $x=\E\mathcal D_\eta$,
\[
 x\le C\left\{\frac{(\log N)^2}{N^2}
                    +\frac1{N\kappa_N}\right\}
                    +\frac C{\kappa_N}\sqrt{x}.
\]
Young's inequality absorbs the last term as $x/2+C/\kappa_N^2$,
giving the residual bound in~\eqref{eq:critical-input-bounds}.
Level crossing gives the first deficit bound; iterating the
$d$-Lipschitz map $u\mapsto\rho_Nu^d$ and applying weighted
Cauchy--Schwarz gives the second. The residual retains the boundary term
$\rho_N-U_1=(I-1)/N$, where $I$ counts empty queues, and satisfies
$0\le\mathcal R_N\le\sum_{k\ge0}d^{-k-1}=1/(d-1)$.
Its small expectation controls the recursion error; it does not force a
deterministic law at the critical center.

We also use queue--phase factorization. For $\ell$ fixed, $x\in\N_0^\ell$ and
$\|f\|_\infty\le1$, solve the finite-state Poisson equation
\[
 -\mathcal A_\ell\psi(x,\cdot)
 =f(x,\cdot)-\sum_{e\in E^\ell}\pi^{\otimes\ell}(e)f(x,e),
 \qquad \pi^{\otimes\ell}\psi(x,\cdot)=0.
\]
Here $\mathcal A_\ell$ generates $\ell$ independent unit-speed phases.
Its inverse on centered functions gives $\|\psi\|_\infty\le C_\ell$.
Queue jumps affecting this test have rate at most $\ell(d\bc+c_*)$.
Stationarity therefore gives, for every $N\ge\ell$, $0<\rho<1$ and $\kappa>0$,
\begin{equation}\label{eq:phase-factorization-common}
 \left|\E f(Q_{1:\ell},J_{1:\ell})-
 \E\sum_e\pi^{\otimes\ell}(e)f(Q_{1:\ell},e)\right|
 \le C_\ell/\kappa .
\end{equation}

The remaining input controls the spread around the empirical mean
$M=N^{-1}\sum_iQ_i$.
\begin{lemma}[Uniform spread around the empirical mean]
\label{lem:boundary-spread}\label{lem:critical-mod-spread}
For constant service there are $t_d>0$ and $C_d<\infty$, depending only
on $d$, such that, for every $N$ and $\rho\in[1/2,1)$,
\begin{equation}\label{eq:boundary-spread}
 \E\left[\frac1N\sum_i e^{t_d|Q_i-M|}\right]\le C_d.
\end{equation}
For independent service environments there are $t>0$, $C<\infty$ and
$\kappa_0<\infty$, depending only on $A,c,d$, such that, for every $N$,
$\rho\in[1/2,1)$ and $\kappa\ge\kappa_0$, the invariant law satisfies
\[
 H_N=\frac1N\sum_i e^{t|Q_i-M|},\qquad \E H_N\le C.
\]
The constant-service assertion requires no switching condition.
\end{lemma}

The proof uses the exponential Lyapunov function
$\Phi=\sum_i(e^{t(Q_i-M)}+e^{-t(Q_i-M)})$. Routing toward shorter
queues gives $\mathcal L^0\Phi\le-a\Phi+bN$ for small fixed $t>0$.
The bounds $|B_\Phi|+|\mathcal L_{\rm q}B_\Phi|\le C\Phi$ show that,
for sufficiently large fixed $\kappa$, the corrected function
$\Phi+B_\Phi/\kappa$ is comparable to $\Phi$ and has the same type of
drift bound by~\eqref{eq:critical-corrected-test}. Bounded concave
truncations then establish integrability and the stated expectations.

\subsection{From the intermediate amplitude to the empirical mean}
We now work under the critical scaling. To verify the family introduced in
Section~\ref{sec:model}, the two sums defining $r(\omega)$ converge locally
uniformly for $\omega>0$, and
\[
 r'(\omega)=-\sum_{j\in\mathbb Z}d^j e^{-\omega d^j}<0,
 \qquad r(d\omega)=r(\omega)-1.
\]
Thus $r$ is a smooth decreasing bijection onto $\mathbb R$; translation by
one yields the inverse relation and bounds in~\eqref{eq:critical-inverse-mean}.
Those bounds also prove that $p_\theta$ in~\eqref{eq:critical-mixing-density}
is integrable at both ends and has a positive finite normalizing constant.

The next estimates compare stationary tests for the averaged and actual
systems. Their errors remain explicit in $N$ and $\kappa$.
\begin{lemma}[Averaging the tests for the critical coordinate]
\label{lem:critical-averaged-tests}
Under $\mu_{N,\rho,\kappa}$, for $N\ge1$, $0<\rho<1$ and $\kappa>0$,
put $X=M_N-m_N$ and $A_k=\sum_i(k-Q_i)_+$, $k\in\N_0$.
For $f\in C_b^2(\mathbb R)$, stationarity yields
\begin{equation}\label{eq:critical-averaged-low-test}
 \big|\E\mathcal L^0[A_kf(X)]\big|
       \le C_f N(k+1)/\kappa .
\end{equation}
For $f$ with bounded $f'$, $f''$ and $f'''$, stationarity gives
\begin{equation}\label{eq:critical-averaged-mean-test}
 \big|\E\mathcal L^0f(X)\big|\le C_f/\kappa .
\end{equation}
The second assertion allows $f$ itself to grow at most linearly.
At $\rho=1-N^{-1}$, $N\ge2$, these estimates imply
\begin{equation}\label{eq:critical-idle-mean}
 \E I=1+O(N/\kappa),\qquad I=\sum_i\one_{\{Q_i=0\}},
\end{equation}
and, for every $h\in C_b^2(\mathbb R)$,
\begin{equation}\label{eq:critical-slow-generator}
 \E[(I-1)h(X)+h'(X)]
       =O_h(N^{-1}+N/\kappa).
\end{equation}
\end{lemma}
\begin{proof}\label{prf:lem:critical-averaged-tests}
Use~\eqref{eq:critical-corrected-test}. Write
$a_f(X)=f(X-N^{-1})-f(X)$ and
$C_r=\sum_i h_{J_i}\one_{\{1\le Q_i\le r\}}$,
with $C_\infty=\sum_i h_{J_i}\one_{\{Q_i>0\}}$. Then
\[
 B_f=a_f(X)C_\infty,\qquad
 B_{A_kf}=A_ka_f(X)C_\infty+f(X-N^{-1})C_k.
\]
A queue jump changes either $C$ by $O(1)$, including at zero,
while $a_f=O_f(N^{-1})$ and its increments are $O_f(N^{-2})$.
Since $A_k\le Nk$ has increments bounded by one,
$B_{A_kf}$ has increments $O_f(k+1)$. The $O(N)$ total queue-jump
rate yields $|\mathcal L_{\rm q}B_{A_kf}|\le C_fN(k+1)$,
proving~\eqref{eq:critical-averaged-low-test} with bounded tests.
Similarly, $|B_f|\le C_f$ and
$|\mathcal L_{\rm q}B_f|\le C_f$, proving
\eqref{eq:critical-averaged-mean-test}.
For linearly growing $f$, Lemma~\ref{lem:finite-recurrence} gives
integrability, and bounded corrected-test increments justify stationary
compensation in~\eqref{eq:critical-corrected-test}.

Take $f(x)=x$. Since
$\mathcal L^0X=\bc(I-1)/N$, this proves
\eqref{eq:critical-idle-mean}.
For $h\in C_b^2$, choose $f$ with $f'=h$.
A second-order Taylor expansion gives
\[
 \frac N{\bc}\mathcal L^0 f(X)
  =(I-1)h(X)
       +\left(1-\frac{1+I}{2N}\right)h'(X)+O_h(N^{-1}).
\]
The remainder is uniform: jumps have size $N^{-1}$ and total rate
at most $2N\bc$. Use~\eqref{eq:critical-idle-mean} and
\eqref{eq:critical-averaged-mean-test} to obtain
\eqref{eq:critical-slow-generator}.
\end{proof}

For the following preparatory lemmas and the proof of
Theorem~\ref{th:critical-random-profile}, assume $\rho_N=1-N^{-1}$ and
$\kappa_N/N\to\infty$, omitting the switching condition for constant
service. Expectations of finite-system quantities are under their invariant law.
For sufficiently large $N$, put
\begin{equation}\label{eq:critical-intermediate-amplitude}
 \begin{gathered}
 k_N=\lfloor\tfrac14\log_dN\rfloor,\qquad
 D_l=N(1-U_l),\qquad I=D_1,\\
 Y_N=D_{k_N}/d^{k_N},\qquad
 H_N=N^{-1}\sum_i e^{t|Q_i-M_N|}.
 \end{gathered}
\end{equation}
Here $k_N\ge1$, and $t$ is from the applicable assertion of
Lemma~\ref{lem:boundary-spread}.
Proposition~\ref{prop:critical-inputs} gives, uniformly in $l\ge1$,
\begin{equation}\label{eq:critical-preliminary-moments}
 \E D_l\le Cd^l,\qquad
 \E D_l^2\le Cd^{2l}(\log N)^2,\qquad
 \E Y_N^2\le C(\log N)^2.
\end{equation}
Set $r(0)=+\infty$, so the mean comparison error is infinite when
$Y_N=0$; this event will have vanishing probability.

The next lemma links the contribution of empty queues to $Y_N$ in the
weighted stationary identity for the mean.
\begin{lemma}[Matching the low levels to the intermediate amplitude]
\label{lem:critical-flux-matching}
For every $f\in C_b^2(\mathbb R)$,
\begin{equation}\label{eq:critical-weighted-matching}
 \E[I f(X_N)]-(d-1)\E[Y_N f(X_N)]\longrightarrow0.
\end{equation}
More precisely, its absolute value is at most
\[
 C_f\left\{d^{-k_N}
       +\frac{d^{k_N}(\log N)^2}{N}+\frac N{\kappa_N}\right\}.
\]
In particular $\E Y_N\to1/(d-1)$.
\end{lemma}
\begin{proof}\label{prf:lem:critical-flux-matching}
Write $A_l=\sum_i(l-Q_i)_+=\sum_{r=1}^lD_r$, a bounded function
for each fixed $N,l$. The homogeneous generator satisfies exactly
\begin{equation}\label{eq:critical-low-drift}
 \bc^{-1}\mathcal L^0 A_l
      =D_{l+1}-I-N\rho_N(1-U_l^d).
\end{equation}
Services at lengths $1,\ldots,l$ increase $A_l$ by one at rate
$\bc(D_{l+1}-I)$; arrivals below $l$ decrease it at rate
$N\bc\rho_N(1-U_l^d)$. Taylor's bound on $[0,1]$ gives
\begin{equation}\label{eq:critical-low-linearization}
 \left|N\rho_N(1-U_l^d)-dD_l\right|
       \le C(D_l+D_l^2)/N.
\end{equation}

Apply~\eqref{eq:critical-averaged-low-test} to $A_lf(X_N)$ and
expand the homogeneous generator of the product. Since $X_N$ jumps
by $\pm N^{-1}$,
\(|\mathcal L^0 f(X_N)|\le C_f(I+1)/N.\)
Minkowski's and Cauchy--Schwarz inequalities, using
\eqref{eq:critical-preliminary-moments}, imply
\(\E A_l^2\le Cd^{2l}(\log N)^2,\qquad \E[A_l(I+1)]\le Cd^l(\log N)^2.\)
The product cross term has expectation at most $C_fd^l/N$:
contributing jumps have rate at most $C(D_{l+1}+D_l)$ and increments
bounded by $1$ and $\|f'\|_\infty/N$.
Together with~\eqref{eq:critical-low-linearization}, this proves
\begin{equation}\label{eq:critical-weighted-recursion}
 \begin{gathered}
 \E[D_{l+1}f(X_N)]
  =d\E[D_lf(X_N)]+\E[I f(X_N)]+R_{l,N},\\
 |R_{l,N}|\le C_f\left\{
       \frac{d^{2l}(\log N)^2}{N}
                      +\frac{N(l+1)}{\kappa_N}\right\}.
 \end{gathered}
\end{equation}
The starting identity $D_1=I$ is exact. Iterate to $k=k_N$ and
divide by $d^k$. The main term is
$(1-d^{-k})\E[I f(X_N)]/(d-1)$; the error is bounded by
\[
 C_f\sum_{l=1}^{k-1}d^{-l-1}
  \left\{\frac{d^{2l}(\log N)^2}{N}
                      +\frac{N(l+1)}{\kappa_N}\right\}
 \le C_f\left\{\frac{d^k(\log N)^2}{N}
                      +\frac N{\kappa_N}\right\}.
\]
Here $\sum_l(l+1)d^{-l}<\infty$; conclude using
$\E I=1+O(N/\kappa_N)$ from~\eqref{eq:critical-idle-mean}.
\end{proof}

We next recover the actual empirical mean from the central tails.
The weighted bounds below allow this relation to enter stationary
balance and will also establish tightness.
\begin{lemma}[An intermediate level determines the profile mean]
\label{lem:critical-shape-mean}
Use the critical scaling and notation fixed above, and abbreviate
$X=X_N$ and $Y=Y_N$. Under the finite-system invariant law,
\begin{equation}\label{eq:critical-Y-moments}
 \E Y\le C,\qquad \E Y^2\le C(\log N)^2,
\end{equation}
and, for every fixed $J<\infty$,
\begin{equation}\label{eq:critical-forward-shape}
 \E\max_{|j|\le J}
   |U_{m_N+j}-e^{-Yd^{j-\theta_N}}|\longrightarrow0.
\end{equation}
There are events $G_N\subset\{Y>0\}$ and deterministic $t_N\downarrow0$
such that
\begin{equation}\label{eq:critical-weighted-mean-match}
 \E[Y\one_{G_N^c}]\longrightarrow0,\qquad
 |X-r(d^{-\theta_N}Y)|\le t_N\quad\hbox{on }G_N.
\end{equation}
For every bounded measurable function $h_R$ with $0\le h_R\le1$ and
$h_R=0$ on $(-\infty,R]$,
\begin{equation}\label{eq:critical-right-weighted-cutoff}
 \limsup_{N\to\infty}\E[Yh_R(X)]\le C_d d^{-R}.
\end{equation}
For every continuous compactly supported function $h$,
\begin{equation}\label{eq:critical-compact-amplitude-match}
 \E[Yh(X)]-\E[d^{\theta_N}\omega(X)h(X)]\longrightarrow0.
\end{equation}
None of these assertions presupposes tightness of $X$.
\end{lemma}
\begin{proof}\label{prf:lem:critical-shape-mean}
Write $\ell=\log N$, $k=k_N$ and $m=m_N$.
Proposition~\ref{prop:critical-inputs} gives
\begin{equation}\label{eq:critical-shape-inputs}
 \E(1-U_l)\le C\min\{1,d^l/N\},\qquad
 \E\mathcal R_N\le C\ell^2/N^2,
\end{equation}
and both moment bounds in~\eqref{eq:critical-Y-moments}.
Lemma~\ref{lem:boundary-spread} gives, for a fixed $t>0$,
\begin{equation}\label{eq:critical-shape-spread-input}
 { H_N=N^{-1}\sum_i e^{t|Q_i-M_N|},\qquad \E H_N\le C.}
\end{equation}

For $s\ge k$, iterate the $d$-Lipschitz map $u\mapsto\rho_Nu^d$ and
use Cauchy--Schwarz with the weights in~\eqref{eq:critical-residual}. This gives
\begin{equation}\label{eq:critical-intermediate-recursion}
 \left|U_s-
 \rho_N^{(d^{s-k}-1)/(d-1)}U_k^{d^{s-k}}\right|^2
 \le C d^{2s-k}\mathcal R_N.
\end{equation}
Indeed the unsquared error is bounded by
$\sum_{l=k}^{s-1}d^{s-1-l}|\rho_NU_l^d-U_{l+1}|$, and the second
squared factor is $\sum_{l=k}^{s-1}d^{2s-1-l}\le Cd^{2s-k}$.

Fix $A>2$, $p>2$, $q>4$, and subsequently choose $B$ sufficiently large.
Set $w_N=\lceil B\log\ell\rceil$ and
\(G_N=\{\ell^{-A}\le Y\le\ell^A,\ H_N\le\ell^p, \mathcal R_N\le\ell^q/N^2\}.\)
On the larger event $\{Y\le\ell^A,\mathcal R_N\le\ell^q/N^2\}$,
including $Y=0$, uniformly for $|j|\le w_N$,
\begin{equation}\label{eq:critical-log-window-shape}
 \begin{gathered}
 |U_{m+j}-e^{-Yd^{j-\theta_N}}|\le e_N,\\
 e_N\le C\left[
 N^{-1/8}\ell^{q/2+B\log d}
 +N^{-3/4}\ell^{2A+B\log d}
 +N^{-1/4}\ell^{B\log d}\right].
 \end{gathered}
\end{equation}
Indeed, $v=Yd^k/N=o(1)$; replacing $(1-v)^{d^{s-k}}$ by
$\exp(-d^{s-k}v)$ contributes an error of $Cd^{s-k}v^2$, and dropping the $\rho_N$
factor contributes an error of $Cd^{s-k}/N$.
Equation~\eqref{eq:critical-intermediate-recursion} bounds the remaining error.
Rounding affects only the constants; all
levels are positive for large $N$. Markov's inequality and
\eqref{eq:critical-Y-moments} make the larger event's complement
negligible. Boundedness of tails then proves
\eqref{eq:critical-forward-shape}.

We next prove the mean assertion on $G_N$. The reference law with tails
$e^{-Yd^{j-\theta_N}}$ puts at least $3/4$ of its mass in
$[-w_N^0,w_N^0]$ for an integer
$w_N^0\le (A/\log d)\log\ell+C_d$.
For large $N$, $w_N^0+1\le w_N$, so
\eqref{eq:critical-log-window-shape} puts at least half of the actual
queues in the same interval about $m$. Consequently
\(H_N\ge\tfrac12e^{t(|X|-w_N^0)_+},\qquad |X|\le K_N:=C_0\log\ell+C_1, \quad C_0=A/\log d+p/t.\)
If $w_N>K_N+1$, then
\begin{align*}
 N^{-1}\sum_i |Q_i-m|\one_{\{|Q_i-m|>w_N\}}
 &\le C(1+K_N)H_N e^{-t(w_N-K_N)/2}.
\end{align*}
Use $|Q_i-m|\le|Q_i-M_N|+K_N$ and absorb the linear factor into
half the exponential.
 With $B>C_0+2(p+2)/t$ and $B\log d>A+2$,
this bound vanishes uniformly on $G_N$. The reference left
first-moment tail is at most $C(w_N+1)\ell^{A-B\log d}$, and its
right tail decays faster than any power of $\ell^{-1}$. Clipped means
differ by at most $(2w_N+1)e_N$, giving the uniform bound
$|X-r(d^{-\theta_N}Y)|\le t_N\to0$.
Taking a decreasing upper envelope of these deterministic errors gives
$t_N\downarrow0$.

For the exceptional sets, Cauchy--Schwarz, Markov's inequality and
\eqref{eq:critical-Y-moments} give the explicit estimate
\begin{equation}\label{eq:critical-shape-bad-sets}
 \E[Y\one_{G_N^c}]
 \le C\left[
 \ell^{-A}+\ell^{2-A}+\ell^{1-p/2}+\ell^{2-q/2}\right]=o(1).
\end{equation}
The terms respectively cover $Y<\ell^{-A}$, $Y>\ell^A$,
$H_N>\ell^p$, and $\mathcal R_N>\ell^q/N^2$, proving
\eqref{eq:critical-weighted-mean-match} without a uniform second moment.

Directly from its convergent tail sums, $r(d\omega)=r(\omega)-1$ and
$|r(\omega)+\log_d\omega|\le C_d$. If $G_N$ occurs and $h_R(X)>0$, then
$r(d^{-\theta_N}Y)>R-t_N$ and hence $Y\le C_d d^{-R}$ for large $N$.
Equation~\eqref{eq:critical-shape-bad-sets} now proves
\eqref{eq:critical-right-weighted-cutoff}.

For the last assertion we first show, for every fixed $K$,
\begin{equation}\label{eq:critical-shape-local-bad}
 \PP(G_N^c,|X|\le K)\longrightarrow0.
\end{equation}
The moment bounds control all exceptional sets except $Y<\ell^{-A}$.
For fixed $J>K$, the spread estimate gives
\(\E[U_{m+J}\one_{\{|X|\le K\}}] \le C e^{-t(J-K)}.\)
On $Y<\ell^{-A}$, the reference tail
$e^{-Yd^{J-\theta_N}}$ tends uniformly to one.
By~\eqref{eq:critical-forward-shape},
\(\limsup_N\PP(Y<\ell^{-A},|X|\le K) \le C e^{-t(J-K)}.\)
Letting $J\to\infty$ proves~\eqref{eq:critical-shape-local-bad}.
If $h$ is supported in $[-K,K]$, then on $G_N\cap\{|X|\le K\}$
the mean approximation and uniform continuity of $\omega$ on
$[-K-1,K+1]$ give $|Y-d^{\theta_N}\omega(X)|=o(1)$ uniformly.
On $G_N^c$, the $Y$ term vanishes by~\eqref{eq:critical-shape-bad-sets}
and the bounded term by~\eqref{eq:critical-shape-local-bad}, proving
\eqref{eq:critical-compact-amplitude-match} without prior tightness.
\end{proof}

\begin{proof}[Proof of Theorem~\ref{th:critical-random-profile}]\label{prf:th:critical-random-profile}
\emph{Step 1: tightness and the mean--amplitude relation.}
We first prove tightness of $X_N$, without assuming uniform
integrability of $Y_N$. From Lemma~\ref{lem:boundary-spread} and
the deficit estimate, for fixed large $R$ and all sufficiently large
$N$,
\begin{equation}\label{eq:critical-left-tightness}
 \PP(X_N\le-R)\le Cd^{-R/2}+Ce^{-tR/2}.
\end{equation}
Take $j=m_N-\lfloor R/2\rfloor$. Then
$\PP(U_j<1/2)\le2\E(1-U_j)\le Cd^{-R/2}$. Otherwise, $X_N\le-R$
places at least half the queues $R/2-1$ above their mean, giving
$H_N\ge\tfrac12e^{t(R/2-1)}$. Markov's inequality proves the bound.

Choose a smooth nondecreasing function $h$ equal to zero on
$(-\infty,0]$ and one on $[1,\infty)$, with
$0\le h'\le K_0$, and let $h_R(x)=h(x-R)$.
Equation~\eqref{eq:critical-right-weighted-cutoff} gives
\begin{equation}\label{eq:critical-upper-weighted-amplitude}
 \E[Y_Nh_R(X_N)]\le Cd^{-R}+o(1)
\end{equation}
for fixed $R$. Equations~\eqref{eq:critical-slow-generator} and
\eqref{eq:critical-weighted-matching} give
$\E h_R(X_N)=\allowbreak\E h'_R(X_N)+\allowbreak
(d-1)\E[Y_Nh_R(X_N)]+\allowbreak o(1)$.
Thus $F_N(R)=\PP(X_N\ge R)$ satisfies
$(1+K_0)F_N(R+1)\le\allowbreak K_0F_N(R)+\allowbreak Cd^{-R}+o(1)$.
Taking upper limits in $N$ and iterating over integers $R\ge0$
gives $\lim_{R\to\infty}\limsup_NF_N(R)=0$; together with
\eqref{eq:critical-left-tightness}, this proves tightness.

The bounded first moments make $Y_N$ tight. To exclude zero, fix
large $R>0$ and then small $\delta>0$. The central-tail approximation gives
\(\PP(Y_N\le\delta,\ U_{m_N+\lceil2R\rceil}<1/2)=o(1),\)
because $e^{-\delta d^{\lceil2R\rceil-\theta_N}}\ge3/4$
when $\delta$ is small enough. If $X_N\le R$ and the empirical
tail in this display is at least $1/2$, then
$H_N\ge\tfrac12e^{tR}$ up to a fixed rounding factor. Consequently
\(\limsup_N\PP(Y_N\le\delta) \le\limsup_N\PP(X_N>R)+Ce^{-tR}.\)
For each fixed $R$ this bound holds for all sufficiently small $\delta$.
Taking $\delta\downarrow0$ and then $R\to\infty$ gives
\(\lim_{\delta\downarrow0}\limsup_N\PP(Y_N\le\delta)=0\).
Thus every subsequential limit of $Y_N$ lies in $(0,\infty)$.
For any $\gamma>0$,
\(\PP(G_N^c)\le\PP(Y_N<\gamma)+\gamma^{-1}\E[Y_N\one_{G_N^c}]\).
Taking $N\to\infty$ and then $\gamma\downarrow0$ gives
$\PP(G_N^c)\to0$. Together with~\eqref{eq:critical-weighted-mean-match},
this yields
\begin{equation}\label{eq:critical-mean-amplitude-relation}
 \PP(Y_N=0)+
 \PP\bigl(Y_N>0,\ |X_N-r(Y_Nd^{-\theta_N})|>\delta\bigr)
       \longrightarrow0,\qquad \delta>0.
\end{equation}

\emph{Step 2: identification of the mixing density.}
Tightness of $X_N$ and continuity of $\omega$ on compact sets give
$Y_N-d^{\theta_N}\omega(X_N)\to0$ in probability. We use~\eqref{eq:critical-compact-amplitude-match} to pass to the limit
in weighted expectations. Convergence in probability alone does not
justify this step.

Let $\mu$ be a weak subsequential limit of $X_N$. For $h\in C_c^2(\mathbb R)$,
combine the slow stationary identity, low-level flux matching, and
\eqref{eq:critical-compact-amplitude-match}. The resulting functions of $X_N$
are bounded and compactly supported, so weak convergence and
$\theta_N\to\theta$ yield
\begin{equation}\label{eq:critical-limit-stationary-identity}
 \int_{\mathbb R}\{h'(x)+(g_\theta(x)-1)h(x)\}\,\mu(\dd x)=0,
 \qquad h\in C_c^2(\mathbb R).
\end{equation}

Distributionally, $\mu'=(g_\theta-1)\mu$. For
$G(x)=\int_0^x(g_\theta(v)-1)\,\dd v$, the measure $e^{-G}\mu$
has zero derivative and hence is a nonnegative constant times
Lebesgue measure. The bounds~\eqref{eq:critical-inverse-mean} make
$e^G$ integrable at both ends; normalization gives
$\mu(\dd x)=p_\theta(x)\,\dd x$. This uniquely identifies every
subsequence and proves~\eqref{eq:critical-mean-limit}.

\emph{Step 3: empirical measures and fixed queues.}
Let $\widetilde\zeta_N=N^{-1}\sum_i\delta_{Q_i-m_N}$.
For each fixed $j$, the central-tail approximation and
\eqref{eq:critical-mean-amplitude-relation} imply
\(U_{m_N+j}-e^{-\omega(X_N)d^j}\longrightarrow0 \quad\hbox{in probability}.\)
With arbitrarily high probability $\omega(X_N)$ stays in a compact
subset of $(0,\infty)$, where the reference mass outside $[-R,R]$
vanishes uniformly as $R\to\infty$.
Adjacent-tail differences control masses inside a fixed window and
endpoint tails control the mass outside. Thus
\begin{equation}\label{eq:critical-queue-profile-TV}
 \E d_{\TV}(\widetilde\zeta_N,\Pi_{X_N})\longrightarrow0.
\end{equation}
Boundedness of total variation justifies the passage from convergence
in probability to convergence of expectations.

To prove empirical queue--phase factorization, put
$V_{q,e}=N^{-1}\sum_i\one_{\{Q_i=q\}}
               (\one_{\{J_i=e\}}-\pi_e)$.
The two-tag case of~\eqref{eq:phase-factorization-common} gives
$\E V_{q,e}^2\le N^{-1}+C/\kappa_N$ uniformly: diagonal terms
contribute at most $N^{-1}$, and off-diagonal terms have zero phase
average and expectation bounded by $C/\kappa_N$. Cauchy--Schwarz
inside a finite centered window and~\eqref{eq:critical-queue-profile-TV}
outside give
$\E d_{\TV}(\widetilde\mu_N,
                  \widetilde\zeta_N\otimes\pi)\to0$,
proving~\eqref{eq:critical-random-empirical-TV}.

Continuity of point masses and their unit sum make
$x\mapsto\Pi_{x}\otimes\pi$ continuous in total variation;
continuous mapping gives the random-measure limit. Exchangeability
compares $\ell$ designated pairs to sampling with replacement from
$\widetilde\mu_N$, with error at most $\ell(\ell-1)/(2N)$.
The product inequality replaces each factor by
$\Pi_{X_N}\otimes\pi$ at vanishing expected error. Weak convergence
of $X_N$ then gives convergence of every mixture point mass. Their
unit sum upgrades this to total variation on the countable state
space, proving~\eqref{eq:critical-tagged-mixture}.
For constant service all environment-correction terms vanish, and
Lemma~\ref{lem:boundary-spread} supplies the required spread bound.
\end{proof}

The double-exponential right tails in Theorem~\ref{th:critical-random-profile}
are conditional on the random mean. The unconditional limiting
queue-length distribution $\bar\Pi_\theta=\int_{\mathbb R}\Pi_xp_\theta(x)\,\dd x$
has an exponential right tail:
\[
 \bar\Pi_\theta([j,\infty))\sim C_\theta e^{-j},\qquad
 C_\theta=\mathcal Z_\theta^{-1}\int_{\mathbb R}e^{-y-\omega(y)}\,\dd y\in(0,\infty),
 \qquad j\to\infty,\ j\in\mathbb Z.
\]
Indeed, $\omega(y+j)=d^{-j}\omega(y)$ and
$V_\theta(y+j)-(y+j)=d^{-j}(V_\theta(y)-y)$. Substituting $x=j+y$ gives
\[
 e^j\bar\Pi_\theta([j,\infty))
 =\mathcal Z_\theta^{-1}\int_{\mathbb R}
       e^{-y-\omega(y)-d^{-j}[V_\theta(y)-y]}\,\dd y.
\]
Dominated convergence applies because $V_\theta(y)-y\ge0$ and
$e^{-y-\omega(y)}$ is integrable by~\eqref{eq:critical-inverse-mean}.
This concerns the limiting marginal: first $N\to\infty$ along the
specified subsequence, then $j\to\infty$.

Since $p_\theta$ vanishes at both endpoints, integrating
$p_\theta'=(g_\theta-1)p_\theta$ gives $\E\omega(X_\theta)=d^{-\theta}/(d-1)$.
Since $\omega(X_\theta)$ is nondegenerate, strict Jensen's inequality yields
\[
 \bar\Pi_\theta([j,\infty))
 =\E e^{-d^j\omega(X_\theta)}
 >\exp\{-d^{j-\theta}/(d-1)\},\qquad j\in\mathbb Z.
\]
Thus every limiting tail probability is strictly larger than its
centered mean-field counterpart.

\section{Diffusion of the mean}\label{sec:diffusion-and-starts}
\subsection{Stationary initialization}\label{sec:critical-diffusion}
The stationary theorem identifies the distribution at each fixed time.
The stationary process will serve as the reference for the later
coupling and initial-layer arguments.
For the process limit, the remaining task is to replace accumulated
idleness by the amplitude and identify the martingale variance. We first
establish uniform integrability; it supplies the time-integrated control
needed for tightness, without higher moments of the centered mean.
\begin{proposition}[Stationary diffusion]\label{th:critical-diffusion}
\label{cor:diffusion-profile-time}
Under the scaling of Theorem~\ref{th:critical-random-profile}, start the
finite system in its invariant law. For every finite $T$,
\begin{equation}\label{eq:critical-functional-limit}
 \mathcal X_N\Longrightarrow X\quad\hbox{in }D([0,T],\mathbb R)
 \hbox{ with the }J_1\hbox{ topology},
\end{equation}
where $X$ is the unique nonexplosive solution of
\eqref{eq:critical-diffusion-sde}, with $X(0)$ of density $p_\theta$
and Brownian motion independent of $X(0)$. The limit is stationary.
Moreover,
\begin{equation}\label{eq:diffusion-time-profile-error}
 \E\int_0^T d_{\TV}(\widetilde\mu_N(t),
                 \Pi_{\mathcal X_N(t)}\otimes\pi)\,\dd t\longrightarrow0,
\end{equation}
and $(\mathcal X_N,\widetilde\mu_N)$ converges jointly to
$(X,\Pi_X\otimes\pi)$ in $J_1$ and the integrated total-variation metric.
For constant service, the queue-only conclusions require no switching
condition.
\end{proposition}

On time $Nt/\bc$, arrivals and services each contribute unit variance
in the limit, giving coefficient $\sqrt2$.

For the rest of this subsection, retain the scaling and lattice-phase
assumptions of Proposition~\ref{th:critical-diffusion}, omitting the switching
condition for constant service. For the dynamic estimates put
\[
 k=k_N=\lfloor\tfrac14\log_dN\rfloor,\quad
 D_l=N(1-U_l),\quad I=D_1,\quad Y_N=D_k/d^k,
 \quad a_k=\frac{1-d^{-k}}{d-1}.
\]
Take $N$ large enough that $k_N\ge2$. Write
$\mathcal I_N(t)=I(Nt/\bc)$ and $\mathcal Y_N(t)=Y_N(Nt/\bc)$
for these observables on the accelerated clock. Other time-dependent
observables below use the same clock.

\begin{lemma}[Stationary input for time-integrated estimates]
\label{lem:diffusion-stationary-input}
The nonnegative variables $Y_N$ are uniformly integrable, and
\begin{equation}\label{eq:diffusion-local-L1}
 \E\left[
  |(d-1)Y_N-g_{\theta_N}(\mathcal X_N(0))|
                  \one_{\{|\mathcal X_N(0)|\le R\}}\right]\to0
 \quad\hbox{for every finite }R.
\end{equation}
Consequently, if $\mathcal Y_N(t)$ denotes the accelerated amplitude,
then the continuous processes
$\int_0^t\mathcal Y_N(s)\,\dd s$ are tight with continuous limits.
More precisely, for every $\varepsilon>0$,
\begin{equation}\label{eq:diffusion-integrated-UI-modulus}
 \lim_{\delta\downarrow0}\limsup_N
 \PP\left(\sup_{\substack{0\le s\le t\le T\\ t-s\le\delta}}
     \int_s^t\mathcal Y_N(u)\,\dd u>\varepsilon\right)=0.
\end{equation}
\end{lemma}
\begin{proof}\label{prf:lem:diffusion-stationary-input}
Theorem~\ref{th:critical-random-profile} and
\eqref{eq:critical-mean-amplitude-relation} imply
$Y_N\Rightarrow d^\theta \omega(X_\theta)$.
Lemma~\ref{lem:critical-flux-matching} gives
$\E Y_N\to1/(d-1)$. The tail bounds for $p_\theta$ make
$g_\theta p_\theta$ integrable and $p_\theta$ vanish at both endpoints.
Integrating $p_\theta'=(g_\theta-1)p_\theta$ therefore gives
$\E g_\theta(X_\theta)=1$, or
$\E[d^\theta\omega(X_\theta)]=1/(d-1)$.
For nonnegative variables, convergence in law and of first moments
gives uniform integrability. On $|\mathcal X_N(0)|\le R$,
\eqref{eq:critical-mean-amplitude-relation} and uniform continuity of
$\omega$ on compact intervals make the localized error in
\eqref{eq:diffusion-local-L1} vanish in probability. Uniform
integrability of $Y_N$ and boundedness of $g_{\theta_N}$ on $[-R,R]$
then give convergence in $L^1$.

For every $K>0$, pathwise,
\[
 \sup_{t-s\le\delta}\int_s^t\mathcal Y_N(u)\,\dd u
 \le K\delta+
       \int_0^T\mathcal Y_N(u)\one_{\{\mathcal Y_N(u)>K\}}\,\dd u.
\]
The last term has expectation $T\E[Y_N\one_{\{Y_N>K\}}]$ by
stationarity. Uniform integrability and then $\delta\downarrow0$
prove~\eqref{eq:diffusion-integrated-UI-modulus}; bounded terminal
expectations complete tightness. The pathwise supremum also covers
stopping-time intervals without any independence assumption.
\end{proof}

\begin{lemma}[Dynamic replacement of the empty-queue count]
\label{lem:diffusion-idle-replacement}
For every finite $T$,
\begin{equation}\label{eq:diffusion-idle-replacement}
 \sup_{0\le t\le T}
 \left|\int_0^t
    \{\mathcal I_N(s)-(d-1)\mathcal Y_N(s)\}\,\dd s\right|
       \longrightarrow0\quad\hbox{in probability}.
\end{equation}
\end{lemma}
\begin{proof}\label{prf:lem:diffusion-idle-replacement}
For $A_l=\sum_i(l-Q_i)_+$, define the nonnegative separable function
\begin{equation}\label{eq:diffusion-low-corrector}
 F_k=\sum_{l=1}^{k-1}d^{-l-1}A_l.
\end{equation}
The exact drift~\eqref{eq:critical-low-drift} telescopes to
\begin{equation}\label{eq:diffusion-low-corrector-drift}
 \bc^{-1}\mathcal L^0F_k=Y_N-a_kI+R_k,
 \qquad
 R_k=\sum_{l=1}^{k-1}d^{-l-1}
                   \{dD_l-N\rho_N(1-U_l^d)\}\ge0.
\end{equation}
The nonnegativity uses $1-u^d\le d(1-u)$ and $\rho_N\le1$.
The estimates in~\eqref{eq:critical-preliminary-moments} give
\begin{equation}\label{eq:diffusion-low-corrector-moments}
 \E F_k\le Ck,
 \qquad
 \E R_k\le C d^k(\log N)^2/N=o(1).
\end{equation}

To account for the service environment, use the corrector
$\widehat F_k=F_k+B_{F_k}/\kappa_N$ from
\eqref{eq:critical-corrected-test}. It satisfies
\begin{equation}\label{eq:diffusion-low-env-bounds}
 |B_{F_k}|\le CN,\qquad
 |\mathcal L_{\rm q}B_{F_k}|\le CN,
 \qquad
 \E\Gamma(\widehat F_k)\le C(1+\kappa_N^{-1}),
\end{equation}
where $\Gamma(F)=\sum_{z'}q(z,z')[F(z')-F(z)]^2$ is the jump
quadratic-variation rate of the full chain.
For the last bound, put
$w_k(q)=\one_{\{1\le q<k\}}\sum_{l=q}^{k-1}d^{-l-1}$.
Departures and arrivals change $F_k$ by $w_k(q)$ and $-w_k(q+1)$,
respectively, both bounded by $Cd^{-q}$. Also
$B_{F_k}=\sum_i h_{J_i}w_k(Q_i)$ is separable, so queue jumps change
only their own summand. Since
\[
 \E\#\{i:Q_i=r\}\le\E D_{r+1}\le Cd^{r+1},
 \qquad
 \sum_{r\ge0}d^{-2r}d^{r+1}<\infty,
\]
the expected queue brackets of $F_k$ and $B_{F_k}$ are bounded.
The phase bracket of $B_{F_k}$ is at most $C\kappa_N$ in expectation;
dividing by $\kappa_N^2$ proves the last bound.
The same increments and bounded server rates give the uniform
bounds in~\eqref{eq:diffusion-low-env-bounds}. For constant service
$B_{F_k}=0$.

Dynkin's formula for $\widehat F_k$, accelerated by $N/\bc$ and
divided by $a_kN$, yields
\begin{equation}\label{eq:diffusion-low-dynkin}
 \begin{split}
 W_N(t)-W_N(0)
  &=\mathcal M_N^F(t)
      +\int_0^t\{\mathcal Y_N(s)/a_k-\mathcal I_N(s)\}\,\dd s
      +E_N^F(t),\\
 W_N(t)&=\widehat F_k(Nt/\bc)/(a_kN).
 \end{split}
\end{equation}
Stationarity and~\eqref{eq:diffusion-low-env-bounds} imply
\[
 \E\sup_{t\le T}|\mathcal M_N^F(t)|^2\le C_T/N,
 \qquad
 \E\operatorname{TV}_{[0,T]}(E_N^F)
    \le C_T\left\{\frac{d^k(\log N)^2}{N}
                                  +\frac N{\kappa_N}\right\}\to0.
\]
Use Doob's inequality for the martingale; the error contains
$R_k/a_k$ and $\mathcal L_{\rm q}B_{F_k}/(a_k\bc\kappa_N)$.

For uniform endpoint control, $F_k\ge0$ gives
$W_N(t)\ge-C/\kappa_N$, whereas stationarity gives
$\E|W_N(t)|\le C(k/N+\kappa_N^{-1})\to0$ at fixed times.
Since $\mathcal I_N\ge0$,~\eqref{eq:diffusion-low-dynkin} yields
\[
 W_N(t)-W_N(s)
 \le\mathcal M_N^F(t)-\mathcal M_N^F(s)
       +a_k^{-1}\int_s^t\mathcal Y_N(u)\,\dd u
       +\operatorname{TV}_{[0,T]}(E_N^F).
\]
On a fixed finite mesh, $W_N$ vanishes in probability. Between mesh
points, this one-sided bound, the vanishing martingale and error, and
\eqref{eq:diffusion-integrated-UI-modulus} control upward excursions.
Let $N\to\infty$ and then the mesh size vanish; the lower bound gives
$\sup_{t\le T}|W_N(t)|\to0$ in probability. Equation~\eqref{eq:diffusion-low-dynkin}
now replaces $\mathcal I_N$ by $\mathcal Y_N/a_k$. Finally,
$a_k^{-1}-(d-1)=O(d^{-k})$ and
$\E\int_0^T\mathcal Y_N(s)\,\dd s\le C_T$ prove the assertion.
\end{proof}

\begin{proof}[Proof of Proposition~\ref{th:critical-diffusion}]\label{prf:th:critical-diffusion}
Apply the environment corrector to the queue-only function
$X=N^{-1}\sum_iQ_i-m_N$. Its explicit form is
\[
 B_X=-\frac1N\sum_i h_{J_i}\one_{\{Q_i>0\}},\qquad
 \widehat X=X+B_X/\kappa_N.
\]
Thus $|B_X|\le C$, $|\mathcal L_{\rm q}B_X|\le C$, and
$\mathcal L^0X=\bc(I-1)/N$.
Integrability, bounded increments and a bounded total rate give a
square-integrable compensated martingale without a second moment of
the initial population. On the accelerated scale,
\begin{equation}\label{eq:diffusion-mean-decomposition}
 \begin{gathered}
 \mathcal X_N(t)=\mathcal X_N(0)+\mathcal M_N(t)
       +\int_0^t\{(d-1)\mathcal Y_N(s)-1\}\,\dd s+E_N(t),\\
 \sup_{t\le T}|E_N(t)|\to0\quad\hbox{in probability}.
 \end{gathered}
\end{equation}
Lemma~\ref{lem:diffusion-idle-replacement} replaces the corrected
drift $I-1$; the other drift error is at most $C_TN/\kappa_N$ and
$|X-\widehat X|\le C/\kappa_N$.
A queue jump changes $\widehat X$ by
$\pm N^{-1}+O((N\kappa_N)^{-1})$, and a phase jump changes it by
$O((N\kappa_N)^{-1})$. Consequently
\begin{equation}\label{eq:diffusion-mean-bracket}
 \langle\mathcal M_N\rangle_t
 =\int_0^t\left\{\rho_N+
   \frac1{N\bc}\sum_i c_{J_i(Ns/\bc)}
                 \one_{\{Q_i(Ns/\bc)>0\}}\right\}\,\dd s
       +O_T(\kappa_N^{-1}).
\end{equation}
Both phase brackets and queue cross terms are uniformly
$O_T(\kappa_N^{-1})$. The stationary autonomous environments have
marginal $\pi^{\otimes N}$, so
\[
 \E\left|N^{-1}\sum_i c_{J_i}-\bc\right|\le C N^{-1/2}.
\]
Removing the empty servers changes this average by at most $c_*I/N$,
whose expectation is $O(N^{-1})$ by
\eqref{eq:critical-idle-mean}. Therefore
\begin{equation}\label{eq:diffusion-bracket-limit}
 \E\sup_{t\le T}|\langle\mathcal M_N\rangle_t-2t|\to0.
\end{equation}
Jumps are bounded by $C/N$, and the accelerated bracket rate is bounded
by a constant independent of the state. Martingale inequalities give uniform fourth
moments on finite intervals and tightness with continuous limits
(see also \citet{Whitt2007}).

The initial laws converge by Theorem~\ref{th:critical-random-profile}.
Lemma~\ref{lem:diffusion-stationary-input} and
\eqref{eq:diffusion-mean-decomposition} therefore imply joint
tightness of $(\mathcal X_N,\mathcal M_N)$ with continuous limits.
In particular, for every $\varepsilon>0$ there is a finite $R$ such
that
\begin{equation}\label{eq:diffusion-compact-containment}
 \limsup_N\PP\left(\sup_{t\le T}|\mathcal X_N(t)|>R\right)
        <\varepsilon.
\end{equation}

Stationarity and~\eqref{eq:diffusion-local-L1} give
\[
 \E\int_0^T
  |(d-1)\mathcal Y_N(s)-g_{\theta_N}(\mathcal X_N(s))|
       \one_{\{|\mathcal X_N(s)|\le R\}}\,\dd s\longrightarrow0.
\]
On $\{\sup_{t\le T}|\mathcal X_N(t)|\le R\}$, the nonnegative
integral inside this expectation bounds the absolute integrated drift
error, uniformly in its upper endpoint. Compact containment
\eqref{eq:diffusion-compact-containment}, followed by $R\to\infty$,
therefore gives global integrated replacement in probability.
No uniform integrability of $g_{\theta_N}(\mathcal X_N(0))$ is needed.

For a subsequential limit $(X,M)$, continuous paths and local uniform
convergence $g_{\theta_N}\to g_\theta$ allow passage in the drift integral.
Equation~\eqref{eq:diffusion-mean-decomposition} therefore gives
\begin{equation}\label{eq:diffusion-limit-integral-equation}
 X(t)=X(0)+M(t)+\int_0^t(g_\theta(X(s))-1)\,\dd s.
\end{equation}

Both $\mathcal M_N$ and $\mathcal M_N^2-\langle\mathcal M_N\rangle$
are martingales in the full filtration, hence admit bounded tests of
past $(\mathcal X_N,\mathcal M_N)$ values. Uniform fourth moments
allow passage against continuous tests of finitely many past values;
a monotone-class argument extends the identities to the limiting
joint filtration. The prelimit bracket need not be adapted to the
smaller filtration. By~\eqref{eq:diffusion-bracket-limit},
$M$ and $M^2-2t$ are martingales there. The continuous martingale
characterization gives $M=\sqrt2 B$ in this filtration; its increments
are independent of the initial sigma-field containing $X(0)$.

The drift $g_\theta-1$ is smooth and locally Lipschitz. For
$V(x)=1+x^2$, its diffusion generator satisfies
\(V''(x)+(g_\theta(x)-1)V'(x) =2+2x(g_\theta(x)-1)\le C(1+x^2).\)
Indeed $2xg_\theta(x)\le0$ for $x<0$, and $g_\theta$ is bounded
for $x\ge0$. Stopping and applying this Lyapunov bound exclude
explosion. Local pathwise uniqueness then identifies every subsequence,
proving~\eqref{eq:critical-functional-limit}. Stationary prelimits and
finite-dimensional convergence give stationarity with density $p_\theta$.

Finally, stationarity makes~\eqref{eq:diffusion-time-profile-error} equal to
$T$ times~\eqref{eq:critical-random-empirical-TV}. At a continuous path, $J_1$ convergence is uniform. The map $x\mapsto\Pi_x\otimes\pi$
is uniformly continuous in total variation on any compact interval
containing the path ranges, so composition is continuous into the
integrated total-variation metric. Continuous mapping and the
vanishing integrated error prove the joint limit.
\end{proof}

\subsection{Transfer from stationary to prepared initial laws}
\label{sec:nonstationary-diffusion}
Preparation means $W_1$-closeness to $\Pi_x$ with the same empirical mean.
To prove Proposition~\ref{th:nonstationary-diffusion}, we compare the given
system with reference systems started from stationary laws conditioned on
a short interval of initial means. Only the initial law is conditioned;
the reference then evolves with the original transition rates and is
generally not stationary. The coupling is nonexpansive after the paired service phases meet;
their earlier discrepancies have uniformly small expected cost.
\begin{proposition}[Transfer from prepared states]
\label{th:nonstationary-diffusion}\label{cor:nonstationary-profile}
Suppose $\rho_N=1-N^{-1}$, $\kappa_N/N\to\infty$ and
$\theta_N\to\theta\in[0,1]$. For an arbitrary initial distribution, put
$Z_N=N^{-1}\sum_iQ_i(0)-m_N$ and $P_N^0=\mu_N^q(0)$. Assume
\begin{equation}\label{eq:nonstationary-initial-preparation}
 Z_N\Longrightarrow\nu,\qquad
 W_1(P_N^0,\Pi_{Z_N})\longrightarrow0\quad\hbox{in probability},
\end{equation}
where $\nu$ is a probability law on $\mathbb R$.
Then $\mathcal X_N$ converges in $D([0,T],\mathbb R)$ with the $J_1$
topology to the diffusion~\eqref{eq:critical-diffusion-sde} started with law
$\nu$, independently of its Brownian motion. For every finite $T$,
\begin{align}
 &\int_0^T W_1(\mu_N^q(t),\Pi_{\mathcal X_N(t)})\,\dd t
       \longrightarrow0\quad\hbox{in probability},
       \label{eq:nonstationary-time-W1}\\
 &\E\int_0^T d_{\TV}(\widetilde\mu_N(t),
          \Pi_{\mathcal X_N(t)}\otimes\pi)\,\dd t\longrightarrow0.
       \label{eq:nonstationary-time-TV}
\end{align}
The pair $(\mathcal X_N,\widetilde\mu_N)$ converges jointly in $J_1$
and integrated total variation. No moment condition on $\nu$, exchangeability, or restriction
on initial phases is required. For constant service, the queue-only
conclusions require no switching condition.
\end{proposition}

\begin{lemma}[Wasserstein preparation and the stationary reference]
\label{lem:nonstationary-W1-input}
For all $x,y\in\mathbb R$,
\begin{equation}\label{eq:nonstationary-profile-W1-isometry}
 W_1(\Pi_{x},\Pi_{y})=|x-y|,
\end{equation}
where $W_1$ uses the distance $|j-j'|$ on the integers.
Under the assumptions and stationary initialization of
Theorem~\ref{th:critical-random-profile}, set
\[
 \mu_N^q=\frac1N\sum_i\delta_{Q_i-m_N},\qquad
 X_N=\frac1N\sum_iQ_i-m_N.
\]
Then
\begin{equation}\label{eq:nonstationary-stationary-W1}
 \E W_1(\mu_N^q,\Pi_{X_N})\longrightarrow0.
\end{equation}
\end{lemma}

\begin{lemma}[A coupling with arbitrary initial service phases]
\label{lem:nonstationary-path-coupling}
Consider two systems with the same $N,\rho,\kappa$ and fixed model
parameters, started from any two queue--phase configurations.
Permute the labels of the second initial configuration so that its
queue lengths are optimally paired with those of the first. There is
a coupling of their subsequent dynamics and a nonnegative variable
$B_N$ such that, for all physical times $s\ge0$,
\begin{equation}\label{eq:nonstationary-all-time-distance}
 \frac1N\sum_i|Q_i^{(1)}(s)-Q_i^{(2)}(s)|
 \le W_1(\mu_N^{q,(1)}(0),\mu_N^{q,(2)}(0))+\frac{B_N}{N},
 \qquad \E B_N\le\frac{CN}{\kappa}.
\end{equation}
The bound is uniform over the two initial configurations; random
initial configurations are allowed by conditioning on them.
No matching condition on the initial phase counts is required.
Moreover, writing
$\Delta_N(s)=N^{-1}\sum_i\one_{\{J_i^{(1)}(s)\ne J_i^{(2)}(s)\}}$,
the accelerated phase discrepancy satisfies
\begin{equation}\label{eq:nonstationary-phase-occupation}
 \E\int_0^T\Delta_N(Nt/\bc)\,\dd t
       \le\frac C{N\kappa},\qquad T<\infty.
\end{equation}
At accelerated time $t$, the empirical joint laws, centered by the same $m_N$, obey
\begin{equation}\label{eq:nonstationary-paired-TV}
 d_{\TV}(\widetilde\mu_N^{(1)}(t),\widetilde\mu_N^{(2)}(t))
 \le\frac1N\sum_i|Q_i^{(1)}(Nt/\bc)-Q_i^{(2)}(Nt/\bc)|
       +\Delta_N(Nt/\bc).
\end{equation}
For constant service one can take $B_N=0$ in the queue-distance
bound. A state-dependent initial permutation preserves the law of
every permutation-invariant reference trajectory observable, even
though it need not preserve the labeled initial law.
\end{lemma}

\begin{lemma}[Reference processes with conditioned stationary initial laws]
\label{lem:nonstationary-conditioned-reference}
Under the scaling and lattice-phase assumptions of
Proposition~\ref{th:nonstationary-diffusion}, fix $T<\infty$.
Let $I$ be a fixed bounded interval with nonempty interior, and start
the finite system in its invariant law conditioned on $X_N\in I$.
Its accelerated mean converges to the solution of
\eqref{eq:critical-diffusion-sde} whose initial law is
$p_\theta(x)\one_I(x)\,\dd x/\int_Ip_\theta$.
The driving Brownian motion remains independent of that initial
value. For this reference process,
\begin{align}
 &\E W_1(\mu_N^q(0),\Pi_{\mathcal X_N(0)})\longrightarrow0,
                                      \label{eq:reference-initial-W1}\\
 &\E\int_0^T
     W_1(\mu_N^q(t),\Pi_{\mathcal X_N(t)})\,\dd t
           \longrightarrow0,          \label{eq:reference-time-W1}\\
 &\E\int_0^T
     d_{\TV}(\widetilde\mu_N(t),(\Pi_{\mathcal X_N(t)}\otimes\pi))\,\dd t
           \longrightarrow0.          \label{eq:reference-time-TV}
\end{align}
The denominator defining the conditioning is positive for all
sufficiently large $N$.
\end{lemma}

The proofs of these three lemmas are given in Supplement, Section~S4.
The coupling bound is uniform over initial phase configurations. We now
give the transfer argument, first for deterministic initial states and then
for arbitrary initial laws.
\begin{proof}[Proof of Proposition~\ref{th:nonstationary-diffusion}]\label{prf:th:nonstationary-diffusion}
We first prove the assertion for deterministic initial states with
$Z_N\to x$ and initial preparation error tending to zero. By the
isometry in Lemma~\ref{lem:nonstationary-W1-input}, this is equivalent
to
\begin{equation}\label{eq:nonstationary-fixed-initial}
 W_1(P_N^0,\Pi_{x})\longrightarrow0.
\end{equation}
Fix $h>0$ and compare with a reference process initialized from the
invariant law conditioned on $X_N\in(x-h,x+h)$. Write $\mathcal X_N'$
and $\mu_N^{q\prime}$ for its accelerated mean and queue measure.
By Lemma~\ref{lem:nonstationary-conditioned-reference}, its initial
empirical measure satisfies
\(W_1(P_N^0,\mu_N^{q\prime}(0)) \le W_1(P_N^0,\Pi_{x})+ W_1(\mu_N^{q\prime}(0),\Pi_{\mathcal X_N'(0)})+h =h+o_{\PP}(1).\)
Apply Lemma~\ref{lem:nonstationary-path-coupling}, ordering queues
with their phases. Its permutation-invariant path-law argument
preserves the reference mean and empirical trajectories. The coupling
gives $D_N\ge0$ such that
\begin{equation}\label{eq:nonstationary-comparison-bound}
 \sup_{t\ge0}|\mathcal X_N(t)-\mathcal X_N'(t)|\le D_N,
 \qquad
 \sup_{t\ge0}W_1(\mu_N^q(t),\mu_N^{q\prime}(t))\le D_N,
 \qquad D_N\le h+o_{\PP}(1).
\end{equation}
Here the additional service error has expectation at most $C/\kappa_N$
and tends to zero. In the constant-capacity case it is zero.

The limiting drift $g_\theta-1$ is decreasing. Under a common Brownian
motion, solutions preserve order and satisfy
\begin{equation}\label{eq:nonstationary-diffusion-contraction}
 \sup_{t\ge0}|X^u(t)-X^v(t)|\le |u-v|.
\end{equation}
Indeed, their difference is absolutely continuous and decreases while
positive; local uniqueness prevents crossing, and nonexplosion extends
the bound to all times. Couple the conditioned limiting initial value
to $x$ using a common independent Brownian motion; the paths then
remain within $h$. Bounded Lipschitz tests, $d_{J_1}\le d_{\rm uniform}$,
Lemma~\ref{lem:nonstationary-conditioned-reference} and
\eqref{eq:nonstationary-comparison-bound}, followed by $h\downarrow0$,
give convergence from $x$.

The triangle inequality and isometry also give
\begin{equation}\label{eq:nonstationary-W1-transfer}
 \int_0^T W_1(\mu_N^q(t),\Pi_{\mathcal X_N(t)})\,\dd t
 \le 2TD_N+
 \int_0^T W_1(\mu_N^{q\prime}(t),
                    \Pi_{\mathcal X_N'(t)})\,\dd t.
\end{equation}
Its right-hand side tends to zero in probability upon first taking
$N\to\infty$ and then $h\downarrow0$.
The family $x\mapsto\Pi_x\otimes\pi$ is $1$-Lipschitz in total variation, because
$d_{\TV}(\xi,\xi')\le W_1(\xi,\xi')$ on $\mathbb Z$.
Lemma~\ref{lem:nonstationary-path-coupling} bounds empirical total
variation by queue and phase discrepancies, the latter with integrated
expectation at most $C/(N\kappa_N)$. Thus
\begin{align}
 &\int_0^T d_{\TV}(\widetilde\mu_N(t),
                             (\Pi_{\mathcal X_N(t)}\otimes\pi))\,\dd t\notag\\
 &\quad\le 2T\min\{1,D_N\}
       +\int_0^T d_{\TV}(\widetilde\mu_N'(t),
                             (\Pi_{\mathcal X_N'(t)}\otimes\pi))\,\dd t\notag\\
 &\qquad+\int_0^T\frac1N\sum_i
          \one_{\{J_i(Nt/\bc)\ne J_i'(Nt/\bc)\}}\,\dd t.
                                      \label{eq:nonstationary-TV-transfer}
\end{align}
This proves convergence of the expected profile error to zero for
every deterministic sequence satisfying
\eqref{eq:nonstationary-fixed-initial}, after $h\downarrow0$.

For a deterministic full state $\sigma=(q,j)$, write
$z_N(\sigma)=N^{-1}\sum_iq_i-m_N$ and\newline
$e_N(\sigma)=W_1(N^{-1}\sum_i\delta_{q_i-m_N},\Pi_{z_N(\sigma)})$.
Let $K_N(\sigma)$ denote its accelerated mean-path law and $K(x)$ the
diffusion law started at $x$. Contraction makes $K$ weakly continuous.
Random starts require uniformity over bounded means and vanishing
preparation errors.
For every compact interval $[-R,R]$ and every deterministic
sequence $\delta_N\downarrow0$,
\begin{equation}\label{eq:nonstationary-uniform-kernel}
 \sup_{\substack{|z_N(\sigma)|\le R\,,\ e_N(\sigma)\le\delta_N}}
 d_{\mathrm{BL}}(K_N(\sigma),K(z_N(\sigma)))\longrightarrow0,
\end{equation}
where $d_{\mathrm{BL}}$ is a bounded Lipschitz metric for weak
convergence on the $J_1$ path space, with a supremum over an empty
set interpreted as zero. If this failed, a counterexample sequence would have
a subsequence with $z_N(\sigma)\to x$, contradicting deterministic-start
convergence and continuity of $K$.
The same argument applied to~\eqref{eq:nonstationary-W1-transfer}
gives, for every $\eta>0$,
\begin{equation}\label{eq:nonstationary-uniform-profile}
 \sup_{\substack{|z_N(\sigma)|\le R\,,\ e_N(\sigma)\le\delta_N}}
 \PP_\sigma\left(\int_0^T
       W_1(\mu_N^q(t),\Pi_{\mathcal X_N(t)})\,\dd t>\eta\right)
       \longrightarrow0.
\end{equation}
It also gives uniform convergence to zero of the expectation in
\eqref{eq:nonstationary-TV-transfer} over these states.

Choose $\delta_N\downarrow0$ so that
$\PP(e_N(Q(0),J(0))>\delta_N)\to0$. Integrate
\eqref{eq:nonstationary-uniform-kernel} over the initial law and use
tightness of $Z_N$ to let $R\to\infty$. The path law is then at
vanishing bounded Lipschitz distance from
$\int K(x)\Law(Z_N)(\dd x)$, which converges by continuity of $K$ to
$\int K(x)\nu(\dd x)$. This is the SDE law with independent initial
value and Brownian motion. The uniqueness and nonexplosion argument
after~\eqref{eq:diffusion-limit-integral-equation} applies conditionally
on each finite initial value, requiring no moment of $\nu$.
Integrating~\eqref{eq:nonstationary-uniform-profile} and the bounded
TV estimate proves the asserted empirical bounds. Since the first limit is continuous
and $x\mapsto\Pi_x\otimes\pi$ is $1$-Lipschitz in total variation, composition is continuous
for the stated topologies. Continuous mapping and the vanishing
expected integrated distance give joint convergence.
\end{proof}

\section{Initial layer}\label{sec:initial-layer}
Theorem~\ref{th:initial-layer-diffusion} removes preparation by combining
two facts: the averaged whole-line dynamics approaches the family $\Pi_x$
at fixed mean, and the finite system follows those dynamics over a growing
physical-time interval. A restart then applies the prepared-state proposition.
The initial assumption is tightness of random measures in $W_1$, so the
argument localizes to compact sets rather than requiring expectations of
their first moments.
\subsection{Relaxation of the whole-line nonlinear dynamics}
\label{subsec:initial-layer-relaxation}
For $\mu\in\mathcal P_1(\mathbb Z)$, put
$v_j(\mu)=\mu([j,\infty))$. The averaged whole-line equations in
physical time are
\begin{equation}\label{eq:initial-layer-whole-line-ODE}
 \dot v_j=\bc\{v_{j-1}^d-v_j^d-v_j+v_{j+1}\},
 \qquad j\in\mathbb Z.
\end{equation}
The whole line has no empty-queue boundary: departures occur at
rate $\bc$ at every integer state. The total upward and downward fluxes of this
limiting equation are both $\bc$, so its mean is conserved. The finite-system
load deficit and empty-queue boundary are retained as error terms in
Supplement S5; their effects vanish on the initial-layer scale, although
the load deficit contributes to the drift on the accelerated scale.

\begin{lemma}[A conservative contractive semiflow]
\label{lem:initial-layer-semiflow}
Equation~\eqref{eq:initial-layer-whole-line-ODE} defines a unique
continuous semiflow $S_s$ on $\mathcal P_1(\mathbb Z)$.
It preserves the mean and satisfies
\begin{equation}\label{eq:initial-layer-semiflow-contraction}
 W_1(S_s\mu,S_s\nu)\le W_1(\mu,\nu),\qquad s\ge0.
\end{equation}
If $\mu\ne\nu$ have the same mean, the inequality is strict for
every $s>0$. The only equilibrium with mean $x$ is
$\Pi_x$.
\end{lemma}

For reference, the tail formula on $\mathbb Z$ reads
\begin{equation}\label{eq:initial-layer-tail-W1}
 W_1(\mu,\nu)=\sum_{j\in\mathbb Z}|v_j(\mu)-v_j(\nu)|,
 \qquad
 \int z\,\mu(\dd z)=\sum_j(v_j(\mu)-v_j^*).
\end{equation}

Here $v_j^*=\one_{\{j\le0\}}$. Quantile coupling and the first-moment
tail formula prove these identities. The construction and strict contraction
of the semiflow are proved in Supplement, Section~S5.
The following relaxation argument stays in the main text: it explains why
an initially unprepared distribution approaches the one-parameter family.
\begin{proposition}[Global relaxation in first Wasserstein distance]
\label{th:initial-layer-relaxation}
For every $\mu\in\mathcal P_1(\mathbb Z)$ of mean $x$,
\begin{equation}\label{eq:initial-layer-global-attraction}
 W_1(S_s\mu,\Pi_{x})\longrightarrow0\qquad(s\to\infty).
\end{equation}
The convergence is uniform on each $W_1$-compact set
$\mathcal C\subset\mathcal P_1(\mathbb Z)$:
\begin{equation}\label{eq:initial-layer-compact-attraction}
 \sup_{\mu\in\mathcal C}
 W_1\left(S_s\mu,\Pi_{\int z\,\mu(\dd z)}\right)
       \longrightarrow0.
\end{equation}
No exponential-moment assumption is imposed in either assertion.
\end{proposition}
\begin{proof}\label{prf:th:initial-layer-relaxation}
First suppose $\mu$ has a finite exponential moment. Contraction
toward the equilibrium with the same mean gives
\begin{equation}\label{eq:initial-layer-first-moment-bound}
 \sup_{s\ge0}\int|z-x|\,S_s\mu(\dd z)
 \le W_1(\mu,\Pi_{x})+
                    \int|z-x|\,\Pi_{x}(\dd z)=:K<\infty.
\end{equation}
Subtracting adjacent equations shows that $p_j=v_j-v_{j+1}$ is the
law of a chain with upward rate
$\bc\beta_j(s)=\bc S_d(v_j(s),v_{j+1}(s))$ and downward rate $\bc$.
Markov's inequality gives, for $j-x\ge R$,
\(\beta_j(s)\le d v_j(s)^{d-1}\le d(K/R)^{d-1}.\)
For $j-x\le-R-1$, it gives instead
\(\beta_j(s)\ge d v_{j+1}(s)^{d-1} \ge d(1-K/R)^{d-1},\qquad R>K.\)
Fix $R$ large enough that these bounds are at most $1/2$ and at
least $(d+1)/2>1$, respectively, uniformly in $s$.

For a sufficiently small $a>0$, consider
$V(j)=e^{a(j-x)}+e^{-a(j-x)}$. Its chain generator is
\begin{align*}
 \bc^{-1}\mathcal A_sV(j)
 &=e^{a(j-x)}\{\beta_j(s)(e^a-1)+e^{-a}-1\}\\
 &\quad+e^{-a(j-x)}\{\beta_j(s)(e^{-a}-1)+e^a-1\}.
\end{align*}
Their derivatives at $a=0$ show that, for sufficiently small $a>0$,
the first coefficient is uniformly negative when $\beta_j\le1/2$,
and the second when $\beta_j\ge(d+1)/2$.
The other exponential is bounded on each half-line, and both are
bounded on the intervening interval. Hence, for $\eta>0,C<\infty$,
\begin{equation}\label{eq:initial-layer-exponential-drift}
 \mathcal A_sV(j)\le-\eta V(j)+C
                    \qquad(j\in\mathbb Z, s\ge0).
\end{equation}
Choose $a$ within the initial exponential-moment exponent. A Poisson
bound with rate $\bc(d+1)$ controls the number of unit jumps, justifying
integration on finite intervals. Gronwall gives
\[
 \int V\,\dd(S_s\mu)
       \le e^{-\eta s}\int V\,\dd\mu+C/\eta.
\]
Thus the orbit is tight with uniformly integrable first moments and
hence relatively compact in $W_1$, by diagonal extraction and tail truncation.
For general $\mu\in\mathcal P_1$, choose finitely supported truncations
$\mu^{(l)}\to\mu$ in $W_1$. Contraction gives
\(\sup_{s\ge0}W_1(S_s\mu^{(l)},S_s\mu)\le W_1(\mu^{(l)},\mu).\)
Their relatively compact orbits approximate the orbit of $\mu$ uniformly
in time, making it totally bounded. Completeness, from the closed
tail-space representation, gives relative compactness.

The decreasing distance $D(s)=W_1(S_s\mu,\Pi_{x})$ tends to
$\ell\ge0$. Along a subsequence $S_{s_n}\mu\to\mu_*$ in $W_1$,
with mean $x$. Continuity and the semiflow property give, for $r\ge0$,
\(W_1(S_r\mu_*,\Pi_{x}) =\lim_n W_1(S_{s_n+r}\mu,\Pi_{x})=\ell.\)
Strict contraction for a distinct same-mean pair contradicts this
for $r>0$. Thus $\mu_*=\Pi_{x}$, $\ell=0$, proving
\eqref{eq:initial-layer-global-attraction}.

For uniformity, define
$D_s(\mu)=W_1(S_s\mu,\Pi_{\int z\,\mu(\dd z)})$.
The contraction and the isometry
$W_1(\Pi_{x},\Pi_{y})=|x-y|$ imply
\(|D_s(\mu)-D_s(\nu)|\le2W_1(\mu,\nu) \qquad\hbox{for all }s\ge0.\)
A finite Wasserstein net of any compact $\mathcal C$, followed
by pointwise convergence at its finitely many centers, proves
\eqref{eq:initial-layer-compact-attraction}.
\end{proof}

\subsection{Finite-system approximation on a growing initial interval}
\begin{lemma}[Approximation during the initial layer]
\label{lem:initial-layer-fluid}
Let $\rho_N=1-N^{-1}$, $m_N=\lfloor\log_dN\rfloor$ and
$\kappa_N/N\to\infty$. Use physical time $s$ and the whole-line
semiflow $S_s$ of Lemma~\ref{lem:initial-layer-semiflow}. Choose any
deterministic sequence
\begin{equation}\label{eq:initial-layer-burnin-scale}
 b_N\longrightarrow\infty,\qquad b_N=o(\log N).
\end{equation}
For every compact set $\mathcal K\subset\mathcal P_1(\mathbb Z)$,
\begin{equation}\label{eq:initial-layer-fluid-compact}
 \sup_{P_N(0)\in\mathcal K}\ \sup_{\text{initial phases}}
 \E\sup_{0\le s\le b_N}W_1(P_N(s),S_sP_N(0))\longrightarrow0,
\end{equation}
where the first supremum is over realizable empirical measures.
If the laws of the random initial profiles $P_N(0)$ are tight in
$\mathcal P_1(\mathbb Z)$, the same supremum over time converges to
zero in probability, without any restriction on the initial phases.
In particular,
\begin{equation}\label{eq:initial-layer-mean-preservation}
 \sup_{0\le s\le b_N}
 \left|\frac1N\sum_iQ_i(s)-\frac1N\sum_iQ_i(0)\right|
 \longrightarrow0\quad\hbox{in probability}.
\end{equation}
For constant service, no switching condition is required.
\end{lemma}

Supplement, Section~S5 gives the quantitative estimate, including the
empty-queue boundary correction and the summed martingale errors.
To see why exponential initial moments are unnecessary, clip centered
queues to $[-R,R]$, keeping their phases. Contractive coupling and
whole-line contraction bound the resulting error by
$2e_R(P_N(0))$, where $e_R(P)=\int(|z|-R)_+P(\dd z)$.
For fixed $R$ and $\eta>0$, the clipped exponential moment is at most
$e^{\eta R}$, so the quantitative error vanishes on $b_N=o(\log N)$.
Then $\sup_{P\in\mathcal K}e_R(P)\to0$ as $R\to\infty$.
Random starts follow by conditioning on the full initial state and
localizing to compact sets of high probability. No lower bound on an
individual service capacity or restriction on the initial phases is used.
On a fixed physical interval $[0,T]$, the same estimate and truncation
argument require only $\kappa_N\to\infty$ under the same initial
$W_1$-tightness assumption.
\subsection{Transfer to the slow stochastic dynamics}
\begin{proof}[Proof of Theorem~\ref{th:initial-layer-diffusion}]\label{prf:th:initial-layer-diffusion}
Uniformly vanishing absolute first-moment tails characterize relative
compactness in $W_1$ on $\mathbb Z$, by diagonal extraction and tail truncation
in~\eqref{eq:initial-layer-tail-W1}. Tightness thus implies
\eqref{eq:initial-layer-tail-condition}. Conversely, for $\varepsilon>0$
choose $\eta_k\downarrow0$ and $R_k\uparrow\infty$ so that
\[
 \sup_N\PP\left(\int |z|\one_{\{|z|>R_k\}}P_N(0)(\dd z)
                                  >\eta_k\right)
             \le\varepsilon 2^{-k},\qquad k\ge1.
\]
The limsup condition gives the bound for large $N$; enlarging $R_k$
handles the finitely many exceptions, whose first moments are almost
surely finite. The set of laws whose $R_k$-tails are at most $\eta_k$ for every $k$
is closed and $W_1$-compact. A union bound gives probability at least
$1-\varepsilon$, proving the converse.

By Lemma~\ref{lem:initial-layer-fluid},
\begin{equation}\label{eq:initial-layer-uniform-fluid-in-probability}
 E_N:=\sup_{s\le b_N}W_1(P_N(s),S_sP_N(0))
                    \longrightarrow0\quad\hbox{in probability}.
\end{equation}
The mean of $S_sP_N(0)$ is $Z_N$, so the difference of means gives
\eqref{eq:initial-layer-main-mean}. Tightness of the initial profiles
and the uniform-on-compact attraction in
Proposition~\ref{th:initial-layer-relaxation} give
\[
 W_1(S_{b_N}P_N(0),\Pi_{Z_N})\longrightarrow0
 \quad\hbox{in probability}.
\]
No relaxation rate is needed here: compact-uniform attraction applies to
every diverging $b_N$, while $b_N=o(\log N)$ is used only in the
finite-system approximation.
The triangle inequality proves~\eqref{eq:initial-layer-preparation}.
Together with~\eqref{eq:initial-layer-main-mean} and
$W_1(\Pi_x,\Pi_y)=|x-y|$, it also gives preparation at the mean at time $b_N$.

Restart the system at the deterministic time $b_N$. The Markov property preserves the
transition law, while~\eqref{eq:initial-layer-main-mean} and
\eqref{eq:initial-layer-preparation} give mean limit $\nu$ and preparation.
Proposition~\ref{th:nonstationary-diffusion} yields the limit of
$\mathcal X_N^{\rm lay}$, including Brownian independence.
Put $\delta_N=\bc b_N/N\to0$. For $t\ge\delta_N$,
$\mathcal X_N(t)=\mathcal X_N^{\rm lay}(t-\delta_N)$.
The continuous limit makes the shifted process's modulus over
$\delta_N$ vanish in probability. For $t<\delta_N$,
\eqref{eq:initial-layer-main-mean} controls
$\mathcal X_N(t)-\mathcal X_N^{\rm lay}(0)$ uniformly; the same modulus
controls $\mathcal X_N^{\rm lay}(t)-\mathcal X_N^{\rm lay}(0)$. Hence
\(\sup_{0\le t\le T}|\mathcal X_N(t)-\mathcal X_N^{\rm lay}(t)| \longrightarrow0\quad\hbox{in probability},\)
proving convergence from time zero.

For the empirical conclusions, Proposition~\ref{th:nonstationary-diffusion} applies after restart at
$b_N$. Changing variables bounds the integrals on $[\delta_N,T]$,
$\delta_N=\bc b_N/N$, by their shifted counterparts on $[0,T]$.
The total-variation integral on $[0,\delta_N]$ is at most $\delta_N$,
also in expectation.

For Wasserstein distance, contraction and mean conservation of the
whole-line semiflow give, for every $s\le b_N$,
\(W_1(P_N(s),\Pi_{\langle z,P_N(s)\rangle}) \le 2E_N+W_1(P_N(0),\Pi_{Z_N}),\)
where $\langle z,P\rangle=\int z\,P(\dd z)$ and $E_N$ is from
\eqref{eq:initial-layer-uniform-fluid-in-probability}.
The second term is a continuous function of $P_N(0)$ and is therefore tight.
Thus the initial-layer supremum is $O_{\PP}(1)$ and its integral
over $\delta_N$ vanishes in probability. Finally, $x\mapsto\Pi_x\otimes\pi$ is
TV-Lipschitz and the mean limit continuous; the vanishing integrated
distance gives joint convergence as in
Proposition~\ref{cor:nonstationary-profile}, without pointwise profile
convergence at zero.
\end{proof}

\section*{Acknowledgments}
The author is also affiliated with the Institute for Information Transmission
Problems of the Russian Academy of Sciences, Moscow, Russia.

The author declares no conflict of interest. No external empirical datasets were used.

\section*{Supplementary material}
\textit{Technical proofs and supporting estimates.} Finite-system recurrence, entropy and spread estimates,
coupling and conditional stationary references, whole-line dynamics,
the quantitative initial-layer approximation, and the equivalent logistic coordinate. See \citet{ASupp}.

\section*{Statement on AI use}

Generative AI tools (DeepSeek) were used for translation into English, language editing, and technical assistance with formatting the manuscript according to the LaTeX template and preparing the reference list. Additionally, DeepSeek was used to assist with code development and preliminary numerical experiments; these preliminary results were used for internal validation and are not reported in this manuscript. The author takes full responsibility for the mathematical statements, references, and final text of the article.

\clearpage
\setcounter{section}{0}
\renewcommand{\thesection}{S\arabic{section}}
\renewcommand{\theHsection}{S\arabic{section}}
\numberwithin{equation}{section}

\pdfbookmark[0]{Supplementary material}{supplement}
\begin{center}
{\Large\bfseries Supplement to ``Critical heavy traffic and diffusion limits for shortest-of-\texorpdfstring{$\routingd$}{d} dynamic routing''\par}
\medskip
Elmira Yu.\ Kalimulina\\
Lomonosov Moscow State University, Moscow, Russia\\
\href{mailto:elmira.yu.k@gmail.com}{\texttt{elmira.yu.k@gmail.com}}
\end{center}
\begin{abstract}
This supplement provides the technical proofs supporting the critical
stationary and diffusion limits in the main paper: finite-system recurrence,
entropy and spread estimates, coupling and conditional stationary references,
the whole-line flow construction, and the quantitative finite-system
approximation used on the initial layer. A final section records an
equivalent coordinate for the limiting diffusion.
\end{abstract}

\section*{Notation and scope}
This supplement gives the technical proofs supporting the main article.
The parameters $d\ge2$, the irreducible generator $A$ on the finite set $E$,
and capacities $(c_e)_{e\in E}$ are fixed. Write $\pi$ for the invariant
phase law, $C=\diag(c_e)$, $\bc=\sum_e\pi_e c_e>0$, $c_*=\max_e c_e$,
$\lambda=\bc\rho$, $\epsilon=1-\rho$ and $a=-\log\rho$.
The $N$ queues have total arrival rate $N\lambda$, shortest-of-$d$
routing with replacement, and independent phase dynamics of generator
$\kappa A$. For a configuration $x$, the arrival rate at queue $i$ is
$\lambda S_d(U_{x_i},U_{x_i+1})$, where
\[
 U_k=N^{-1}\sum_i\one_{\{Q_i\ge k\}},\qquad
 S_d(u,v)=\sum_{r=0}^{d-1}u^{d-1-r}v^r,\qquad
 z_k=\rho^{(d^k-1)/(d-1)}.
\]
The centered phase corrector solves $-Ah=c-\bc\one$ and $\pi^\top h=0$.
The finite invariant law is $\mu_{N,\rho,\kappa}$. The sequence $z_k$
is the classical constant-service mean-field reference used in the entropy
calculation.
Total variation is $d_{\TV}(P,Q)=\sup_B|P(B)-Q(B)|$, and $W_1$
uses integer distance $|j-j'|$.

The critical centering is $m_N=\lfloor\log_dN\rfloor$ with
$\theta_N=\log_dN-m_N$.
Write $M_N=N^{-1}\sum_iQ_i$, $X_N=M_N-m_N$
and $I=\sum_i\one_{\{Q_i=0\}}$.
Physical time $s$ and accelerated time $t$ are related by $s=Nt/\bc$.
Thus $P_N(s)=N^{-1}\sum_i\delta_{Q_i(s)-m_N}$,
$\mu_N^q(t)=P_N(Nt/\bc)$, and $\widetilde\mu_N(t)$ is the corresponding
centered empirical queue--phase measure. Also
$\mathcal X_N(t)=N^{-1}\sum_iQ_i(Nt/\bc)-m_N$.

For $\omega>0$, let $\xi_\omega([j,\infty))=e^{-\omega d^j}$.
Its mean $r(\omega)$ is a strictly decreasing bijection onto $\mathbb R$;
write $\omega(x)=r^{-1}(x)$ and $\Pi_x=\xi_{\omega(x)}$.
The critical drift uses $g_\theta(x)=(d-1)d^\theta\omega(x)$, and
\[
 p_\theta(x)=Z_\theta^{-1}
       \exp\left\{-x-\int_x^\infty g_\theta(v)\,\dd v\right\}.
\]
The whole-line flow $S_s$ has scalar tails $v_j(\mu)=\mu([j,\infty))$
and satisfies $\dot v_j=\bc(v_{j-1}^d-v_j^d-v_j+v_{j+1})$.
In the finite approximation, $V_j^N(s)=P_N(s)([j,\infty))$.
All references to results in the main article are explicitly identified;
local equation numbers refer to this supplement.

\section{Stationary foundations}\label{supp:stationary-foundations}
\subsection{Recurrence of the finite system}
The first proof applies to every $N\ge1$, $0<\rho<1$ and $\kappa>0$.
For its uniform first-moment bound, $\kappa$ is bounded below by a fixed
positive constant.

\begin{proof}[Proof of Lemma~3.1 of the main paper]\label{supp:prf:lem:finite-recurrence}
Fix $\kappa_*>0$ and $\kappa\ge\kappa_*$; for existence at any given
$\kappa>0$, take $\kappa_*=\kappa$. The total jump rate is bounded for
fixed $N,\kappa$, so the chain is nonexplosive.
Writing $\widehat q_k=U_k(x)$, (2) of the main paper gives
$\sum_i b_i=N\lambda$, $b_i\le d\lambda$, and
\[
 \frac1N\sum_i x_i b_i(x)
   =\lambda\sum_{k\ge1}\widehat q_k^d
   \le\frac\lambda N\sum_i x_i.
\]
Put $H=\|h\|_\infty$, $r_e=2(H-h_e)/\kappa\in[0,4H/\kappa]$,
and $V(x,e)=\sum_i(x_i^2+r_{e_i}x_i)$. The Poisson equation gives
\[
 \mathcal LV
 =2\sum_i x_i(b_i-\bc)+\sum_i b_i(1+r_{e_i})
   +\sum_i c_{e_i}\one_{\{x_i>0\}}(1-r_{e_i}).
\]
Indeed, the phase contribution $-2\sum_i x_i(Ah)_{e_i}$ cancels the
phase-dependent part of the linear service drift. The preceding identities
and $0\le r_e\le4H/\kappa_*$ imply
\[
 \mathcal LV\le-2\bc\epsilon\sum_i x_i+NC_0,
 \qquad C_0=d\bc(1+4H/\kappa_*)+c_*.
\]
Cap each queue at $K$, retaining the same sampling rule and rejecting an
arrival routed to a full queue, without rerouting. The rejected increments
$2K+1+r_{e_i}$ of $V$ are nonnegative, so the same drift bound holds.
The capped chain is finite and irreducible. In its invariant law,
$\E\mathcal L_KV=0$ gives $\E\sum_iQ_i\le NC_0/(2\bc\epsilon)$,
uniformly in $K$, hence tightness for fixed $N$. For each fixed state,
there are finitely many incoming neighbors, and all relevant rates equal
the uncapped rates for sufficiently large $K$. Coordinate balance equations
therefore pass to any subsequential limit. Uniformization using the bounded
total jump rate shows that this limiting probability law is invariant;
Fatou preserves the moment bound.

The uncapped chain is irreducible: working phases without arrivals empty
it, and sampling any prescribed server $d$ times permits an arrival there.
Thus its invariant law is unique. Permutation symmetry of the generator
then gives exchangeability, and the moment bound yields
$\E Q_i\le C_0/(2\bc\epsilon)$. Autonomy of the phases gives their
stationary marginal $\pi^{\otimes N}$; queues and phases need not be
independent. The stated first-moment constant is $B=C_0/(2\bc)$.
\end{proof}

\section{Entropy bounds for the critical recursion}\label{supp:entropy}
\subsection{A weighted recursion residual}
We prove the finite-system estimate underlying Proposition~3.2 of the main paper.
The entropy calculation is convenient at general load; its critical
specialization is recorded immediately after the proof.

\begin{proposition}[Weighted recursion residual]\label{prop:supp-entropy}
Suppose $N\ge2$, $\rho\ge1/2$, $\kappa\ge1$ and
$b=\epsilon+\kappa^{-1}\le1/2$. Under the finite-system invariant law, put
\begin{equation}\label{eq:supp-residual-definition}
 \mathcal R_N=\sum_{k\ge0}d^{-k-1}(\rho U_k^d-U_{k+1})^2,
 \qquad
 \mathcal E_{N,\epsilon,\kappa}
 =\frac{b(1+\log(1/b))\log N}{N}
   +\frac{\log N}{N^3}+\frac{b+N^{-1}}{\kappa}.
\end{equation}
There is a constant $C$, depending only on $A,c,d$, such that
\begin{equation}\label{eq:supp-residual-bound}
 \E\mathcal R_N\le C\{\mathcal E_{N,\epsilon,\kappa}+N^{-6}\},
 \qquad
 \E(1-U_l)\le C\min\{1,b d^l\}\quad(l\ge0).
\end{equation}
For constant service, take $b=\epsilon$ and omit the last term in
$\mathcal E_{N,\epsilon,\kappa}$; no switching parameter is required.
\end{proposition}

\begin{proof}\label{prf:prop:supp-entropy}
For fixed $N,\rho,\kappa$, the total jump rate is bounded by
$N(\lambda+c_*+\kappa\max_e(-a_{ee}))$. Thus an integrable test function with bounded jump increments has an
integrable compensated jump sum and
stationary drift $\E\mathcal LF=0$: its expected absolute jump
variation on $[0,t]$ is bounded by $t$ times this rate times
$\sup|\Delta F|$. This applies below without a uniform bound in $N,\kappa$.

We construct a nonnegative entropy dissipation that controls the errors
in the stationary tail recursion. Stationarity then bounds this dissipation
by a finite-system jump remainder and a service-environment correction.
Squared tail deficits control the correction without an additional
logarithmic factor.

\emph{Stationary flux and a deficit bound.}
By Lemma~3.1 of the main paper, the stop-loss test
$N^{-1}\sum_i(Q_i-k)_+$ is integrable and has bounded increments.
Its balance gives, for $k\ge0$,
\begin{equation}\label{eq:entropy-stationary-flux}
 \lambda\E U_k^d
 =\E\frac1N\sum_i c_{J_i}\one_{\{Q_i\ge k+1\}}.
\end{equation}
Apply stationarity also to
$B_k=N^{-1}\sum_i h_{J_i}\one_{\{Q_i\ge k+1\}}$.
Its phase drift is
$-\kappa N^{-1}\sum_i(c_{J_i}-\bc)\one_{\{Q_i\ge k+1\}}$.
Only level-$k+1$ crossings contribute to its queue drift, bounded by
$\|h\|_\infty(d\bc+c_*)$. Hence, uniformly in $k,N,\rho$,
\begin{equation}\label{eq:entropy-approximate-recursion}
 \left|\E U_{k+1}-\rho\E U_k^d\right|\le C/\kappa .
\end{equation}
For constant service, $\E U_{k+1}=\rho\E U_k^d$ exactly.
Jensen's inequality, or $1-x^d\le d(1-x)$, and induction now give
\begin{equation}\label{eq:entropy-deficit}
 \E(1-U_k)\le C\min\{1,b d^k\},\qquad k\ge0.
\end{equation}
Indeed, if $v_k=\E(1-U_k)$, then $v_0=0$ and
$v_{k+1}\le d v_k+\epsilon+C/\kappa$.

\emph{A smoothed entropy and its deterministic drift.}
Put $\eta=N^{-3}$ and, for $0\le u\le1$, define
\[
 \phi_\eta(u)=(u+\eta)\log(u+\eta)-\eta\log\eta-u,
 \qquad
 V_\eta(U)=\sum_{k\ge1}d^{-k}
                 \{\phi_\eta(U_k)-U_k\log z_k\}.
\]
Each summand vanishes at $U_k=0$. The entropy is bounded below and
integrable since $|V_\eta|\le C+Ca\sum_{k\ge1}U_k$.
Its derivatives $G_k=d^{-k}\log((U_k+\eta)/z_k)$ satisfy
$|G_k|\le C\log N$ uniformly; set $G_0=0$.

Let
\[
 \begin{gathered}
 \mathcal J_k=\rho U_k^d-U_{k+1},\quad
 A_k^\eta=\rho(U_k+\eta)^d,\quad B_k^\eta=U_{k+1}+\eta,\\
 \mathcal D_\eta(U)=
 \sum_{k\ge0}d^{-k-1}(A_k^\eta-B_k^\eta)\log\frac{A_k^\eta}{B_k^\eta} .
 \end{gathered}
\]
Every summand of $\mathcal D_\eta$ is nonnegative.
Let $\mathcal L^0$ denote the generator obtained by replacing every service
capacity by $\bc$. Then
$\mathcal L^0 U_k=\bc(\mathcal J_{k-1}-\mathcal J_k)$ for $k\ge1$.
Here $\mathcal J_k$ is the residual in the equilibrium recursion
$U_{k+1}=\rho U_k^d$.
Using $z_{k+1}=\rho z_k^d$, summation by parts gives
\begin{align}
 \sum_{k\ge1}G_k(\mathcal J_{k-1}-\mathcal J_k)
 &=-\sum_{k\ge0}d^{-k-1}\mathcal J_k\log\frac{A_k^\eta}{B_k^\eta}
       +\mathcal J_0\log(1+\eta) \notag\\
 &=-\mathcal D_\eta(U)+O(\eta\log N).
 \label{eq:entropy-deterministic-drift}
\end{align}
The error is uniform because $|A_k^\eta-B_k^\eta-\mathcal J_k|\le C\eta$,
$|\log(A_k^\eta/B_k^\eta)|\le C\log N$ and $|\mathcal J_0|\le1$.
Only finitely many $\mathcal J_k$ are nonzero in each configuration; the remaining
weighted sums converge absolutely, justifying summation by parts.

\emph{The stationary jump remainder.}
Put $\Delta=1/N$ and let
\[
 \beta_k=N\lambda(U_{k-1}^d-U_k^d),\qquad
 \delta_k=\sum_i c_{J_i}\one_{\{Q_i=k\}},\qquad
 \delta_k^0=N\bc(U_k-U_{k+1})
\]
be, respectively, the upward rate of $U_k$, its actual downward rate,
and its homogeneous downward rate.
For the convex function $\phi_\eta$ set
\begin{align*}
 r_+(u)&=\phi_\eta(u+\Delta)-\phi_\eta(u)
                         -\Delta\log(u+\eta),\\
 r_-(u)&=\phi_\eta(u-\Delta)-\phi_\eta(u)
                         +\Delta\log(u+\eta).
\end{align*}
Only arguments for which the associated jump is possible are used.
Both remainders are nonnegative. The resulting generator identity is
\[
 \mathcal L^0 V_\eta
 =\bc\sum_{k\ge1}G_k(\mathcal J_{k-1}-\mathcal J_k)+\mathcal R_\eta^0,
 \quad
 \mathcal R_\eta^0
 =\sum_{k\ge1}d^{-k}
           \{\beta_k r_+(U_k)+\delta_k^0r_-(U_k)\}.
\]
Applying the full generator to $\one_{\{U_k\ge(j+1)\Delta\}}$ gives
exact level crossing even though $U_k$ need not be Markov:
\[
\E[\beta_k\one_{\{U_k=j\Delta\}}]
 =\E[\delta_k\one_{\{U_k=(j+1)\Delta\}}],\qquad 0\le j<N.
\]
Consequently
$\E[\beta_k r_+(U_k)]
 =\E[\delta_k r_+(U_k-\Delta)]$.
At a positive lattice point $u$ there is the exact identity
\[
 r_+(u-\Delta)+r_-(u)
   =\Delta\log\frac{u+\eta}{u-\Delta+\eta}
   \le 3\log N\,\frac{\Delta^2}{u}.
\]
For $u=\Delta$ this follows from $\eta=N^{-3}$; for $u\ge2\Delta$
use $\log(1+x)\le x$ and $u-\Delta\ge u/2$.
Since both $\delta_k$ and $\delta_k^0$ are at most
$N\max(c_*,\bc)(U_k-U_{k+1})$, Tonelli's theorem yields
\[
 \E\mathcal R_\eta^0
 \le \frac{C\log N}{N}\sum_{k\ge1}d^{-k}
       \E\left[\frac{U_k-U_{k+1}}{U_k};\,U_k>0\right].
\]
On $\{U_k\ge1/2\}$ the ratio is at most $2(1-U_{k+1})$;
on $\{0<U_k<1/2\}$ it is at most one and this event has probability
at most $2\E(1-U_k)$.
By~\eqref{eq:entropy-deficit},
\begin{equation}\label{eq:entropy-jump-bound}
 \E\mathcal R_\eta^0
 \le\frac{C\log N}{N}\sum_{k\ge1}d^{-k}\min\{1,b d^{k+1}\}
 \le \frac{Cb(1+\log(1/b))\log N}{N}.
\end{equation}

\emph{Pathwise control of squared deficits.}
For positive $x,y$,
$(x-y)\log(x/y)\ge2(x-y)^2/(x+y)$.
Since $A_k^\eta+B_k^\eta\le2^d+2$ and
$|A_k^\eta-B_k^\eta-\mathcal J_k|\le C\eta$, we have, for every configuration,
\begin{equation}\label{eq:entropy-pathwise-flux}
 \sum_{k\ge0}d^{-k-1}\mathcal J_k^2
 \le C\mathcal D_\eta(U)+C\eta^2.
\end{equation}
Writing $e_k=U_k-z_k$, the recursion for $z_k$ gives
$|e_{k+1}|\le d|e_k|+|\mathcal J_k|$ with $e_0=0$.
Cauchy--Schwarz with geometric weights therefore yields
\begin{equation}\label{eq:entropy-pathwise-profile}
 |e_m|\le d^{m-1}\sum_{k=0}^{m-1}d^{-k}|\mathcal J_k|,
 \qquad
 |e_m|^2\le Cd^{2m}\{\mathcal D_\eta(U)+\eta^2\}.
\end{equation}
Also $1-z_k\le C\epsilon d^k$. Put
\[
 s=\epsilon+\sqrt{\mathcal D_\eta(U)}+\eta,
 \qquad W_2(U)=\sum_{k\ge1}d^{-k}(1-U_{k+1})^2.
\]
It follows pointwise that $1-U_k\le\min\{1,Cd^k s\}$ and hence
\begin{equation}\label{eq:entropy-squared-deficit}
 W_2(U)\le C\{\epsilon+\sqrt{\mathcal D_\eta(U)}+\eta\}.
\end{equation}
For $0<s<1$, split at $K=\lfloor\log_d(1/s)\rfloor$:
the two parts are bounded by $Cs^2d^K\le Cs$ and $Cd^{-K}\le Cds$.
For $s\ge1$, $W_2\le\sum_{k\ge1}d^{-k}\le Cs$.
No concentration of $U$ is assumed.

\emph{Correcting for the service environment.}
For a queue configuration $x$ let
\[
 S_i(x)=\one_{\{x_i>0\}}\{V_\eta(x-e_i)-V_\eta(x)\},
 \qquad
 B(x,e)=\sum_i h_{e_i}S_i(x).
\]
Here $V_\eta(x)$ is the entropy of its empirical tails.
The gradient bound gives $|S_i|\le C\log N/N$ and $|B|\le C\log N$.
Writing $\mathcal L=\mathcal L_{\rm q}+\kappa\mathcal A$ with
$\mathcal A$ acting on all phases, $-Ah=c-\bc\one$ gives
\[
 \mathcal L\left(V_\eta+\frac{B}{\kappa}\right)
 =\mathcal L^0V_\eta+\frac{\mathcal L_{\rm q}B}{\kappa}.
\]
It remains to estimate $\E|\mathcal L_{\rm q}B|$.
We separate the contribution of the server whose queue changes from the
contribution of the other servers at the crossed level. Keeping the level
weights in both contributions will yield the squared-deficit bound
in~\eqref{eq:entropy-corrector-stationary}.
For the remaining queue drift, put
\(F_k(u)=d^{-k}\{\phi_\eta(u-\Delta)-\phi_\eta(u) +\Delta\log z_k\},\qquad u\ge\Delta.\)
Then $S_i=F_{x_i}(U_{x_i})$ for $x_i>0$. Put
$p_k=U_k-U_{k+1}$ for $k\ge0$ and set $F_0=0$.
For every $u\in\{\Delta,2\Delta,\ldots,1\}$,
\begin{equation}\label{eq:entropy-corrector-size}
 |F_k(u)|\le\Delta\left\{d^{-k}\left(-\log u+\frac{\Delta}{u}+\eta\right)+\frac{a}{d-1}\right\}.
\end{equation}
Indeed,
\[
 \Delta^{-1}\{\phi_\eta(u)-\phi_\eta(u-\Delta)\}
 =\int_0^1\log(u-t\Delta+\eta)\,\dd t
 \in[\log u-\Delta/u,\eta],
\]
where the lower bound follows from
$\int_0^1-\log(1-rt)\,\dd t\le r$ for $0\le r\le1$.
The integral is convex in $r$ and has endpoint values zero and one,
so it lies below its chord. The upper bound uses $u\le1$ and
$\log(1+\eta)\le\eta$; also $d^{-k}|\log z_k|\le a/(d-1)$.
The lower-bound singularity at $u=\Delta,t=1,\eta=0$ is integrable.
For lattice points $\Delta\le u\le1-\Delta$,
\begin{equation}\label{eq:entropy-corrector-increment}
 |F_k(u+\Delta)-F_k(u)|\le C d^{-k}\frac{\Delta^2}{u}.
\end{equation}
Indeed the exact second difference is
\[
 \phi_\eta(u+\Delta)-2\phi_\eta(u)+\phi_\eta(u-\Delta)
 =\int_{-\Delta}^{\Delta}\frac{\Delta-|s|}{u+s+\eta}\,\dd s.
\]
For $u\ge2\Delta$ it is at most $2\Delta^2/u$.
At $u=\Delta$ it is maximized at $\eta=0$, where it equals
$2\Delta\log2$.

Retain the level weights in the absolute jump variation of $B$.
At a crossing of level $k$, only $U_k$ changes; $U_{k-1}$ is unchanged.
The contribution of the jumping server changes at an arrival across
level $k$ and at a departure by
\[
 h_e\{F_k(U_k+\Delta)-F_{k-1}(U_{k-1})\},\qquad
 h_e\{F_{k-1}(U_{k-1})-F_k(U_k)\},
\]
respectively. First bound the phase coefficient by $\|h\|_\infty$.
For $f_k=|F_k|$, level crossing gives
\begin{equation}\label{eq:entropy-corrector-level-crossing}
 \E[\beta_k f_k(U_k+\Delta)]
       =\E[\delta_k f_k(U_k)].
\end{equation}
Zero-rate products are zero, and $f_k$ is evaluated only at positive
allowed points. Bound the new birth and old departure terms by
$\delta_k\le Nc_*p_k$ and~\eqref{eq:entropy-corrector-size}.
For the old birth term use $\beta_k\le Nd\lambda p_{k-1}$ and shift
indices; $F_0=0$ removes $k=1$. For the new departure term ($k\ge2$),
$U_{k-1}\ge U_k>0$, both $-\log u$ and $\Delta/u$ decrease, and
$d^{-(k-1)}=d\cdot d^{-k}$. Thus the four contributions total at most
\[
 C\left\{a+\Delta+\eta+\sum_{k\ge1}d^{-k}
             \E[p_k(-\log U_k);\ U_k>0]\right\}.
\]
The $a$ term uses total queue rate $N(\lambda+c_*)$, and the
$\Delta$ term uses $p_k/U_k\le1$ and summable weights.
For $U_k\ge1/2$, $-\log U_k\le2(1-U_k)$ and $p_k\le1-U_{k+1}$;
for $0<U_k<1/2$, $-U_k\log U_k\le1/e$ and $1-U_{k+1}\ge1/2$.
Hence $p_k(-\log U_k)\le C(1-U_{k+1})^2$ in both cases, bounding
the jumping-server contribution by $C\{a+\Delta+\eta+\E W_2(U)\}$.

Among other servers, only those of length $k$ contribute when $U_k$
jumps. At a birth there are $Np_k$; multiplying
\eqref{eq:entropy-corrector-increment} by this count and
$\beta_k\le Nd\lambda p_{k-1}$ bounds their absolute jump variation
rate by $Cd^{-k}p_{k-1}p_k/U_k$, interpreted as zero at $U_k=0$.
At a departure there are $Np_k-1$ other servers. If $U_k=\Delta$,
there are none. Otherwise the corresponding bound, after multiplication
by $\delta_k\le Nc_*p_k$, is
\[
 C d^{-k}p_k\frac{p_k-\Delta}{U_k-\Delta}
 \le C d^{-k}\frac{p_k^2}{U_k},
\]
where the inequality uses $p_k\le U_k$; the expression is used only
when a departure is possible, so $p_k\ge\Delta$.
For $U_k\ge1/2$, both displayed ratios are at most
$2(1-U_{k+1})^2$, since
$p_{k-1}\le1-U_k$ and $p_k\le1-U_{k+1}$.
For $0<U_k<1/2$, both are at most one and hence at most
$4(1-U_{k+1})^2$. Summing bounds the other-server contribution by
$C\E W_2(U)$.

Tonelli, the jumping-server bound,
\eqref{eq:entropy-squared-deficit} and $a\le2\epsilon$ yield
\begin{equation}\label{eq:entropy-corrector-stationary}
 \begin{aligned}
 \E|\mathcal L_{\rm q}B|
 &\le C\{a+\Delta+\eta+\E W_2(U)\}\\
 &\le C\{\epsilon+\Delta+\eta+
                    \E\sqrt{\mathcal D_\eta(U)}\}.
 \end{aligned}
\end{equation}
No queue--phase independence was used. The corrected entropy is
integrable with bounded increments, so its expected drift is zero.
Set $x=\E\mathcal D_\eta(U)<\infty$, since pointwise
$\mathcal D_\eta\le C\log N$. Equations
\eqref{eq:entropy-deterministic-drift}, \eqref{eq:entropy-jump-bound},
\eqref{eq:entropy-corrector-stationary} and Jensen's inequality give
\[
 x\le C\left\{\frac{b(1+\log(1/b))\log N}{N}
       +\eta\log N+\frac{\epsilon+\Delta+\eta}{\kappa}\right\}
       +\frac C\kappa\sqrt{x}.
\]
Young's inequality bounds the last term by $x/2+C/\kappa^2$.
Since $\eta\le\Delta$, $\epsilon/\kappa+\kappa^{-2}=b/\kappa$
and $\eta=N^{-3}$, we obtain
\begin{equation}\label{eq:entropy-dissipation-bound}
 \E\mathcal D_\eta(U)\le C\mathcal E_{N,\epsilon,\kappa}.
\end{equation}
Taking expectations in~\eqref{eq:entropy-pathwise-flux} proves
\eqref{eq:supp-residual-bound}, since $\eta^2=N^{-6}$.
For constant service, omit the correction and use the exact stationary
recursion, giving the asserted bound without the last term.
\end{proof}

\subsection{Critical specialization and pathwise propagation}
For $N\ge4$, $\rho_N=1-N^{-1}$ and $\kappa_N\ge N$,
we have $N^{-1}\le b\le2N^{-1}\le1/2$.
The preceding proposition therefore gives
\begin{equation}\label{eq:supp-critical-residual}
 \E\mathcal R_N\le C\frac{(\log N)^2}{N^2},\qquad
 \E D_l\le C\min\{N,d^l\},\quad D_l=N(1-U_l).
\end{equation}
The same conclusions hold for constant service. These bounds concern the
recursion residual, not concentration about a deterministic law at criticality.

For $s\ge k\ge0$, iterate $\varphi(u)=\rho_Nu^d$. Since $\varphi$
is $d$-Lipschitz on $[0,1]$, Cauchy--Schwarz with weights $d^{-l-1}$ gives
\begin{align}
 \left|U_s-\rho_N^{(d^{s-k}-1)/(d-1)}U_k^{d^{s-k}}\right|^2
 &\le\left(\sum_{l=k}^{s-1}d^{2s-l-1}\right)\mathcal R_N\notag\\
 &=\frac{d^{2s-k}(1-d^{-(s-k)})}{d-1}\,\mathcal R_N.
 \label{eq:supp-residual-propagation}
\end{align}
The case $s=k$ gives zero on both sides. The level $k=0$ is included:
$\mathcal J_0=(I-1)/N$. Also
$0\le\mathcal R_N\le\sum_{k\ge0}d^{-k-1}=1/(d-1)$, and the sum
is finite for each finite queue configuration; hence no moment assumption
is needed to define it. Taking $k=0$ and using
$1-z_l\le Cd^l/N$ yields
\[
 \E D_l^2\le2N^2\E|U_l-z_l|^2+2N^2(1-z_l)^2
       \le Cd^{2l}(\log N)^2.
\]
This proves the residual and deficit estimates of
Proposition~3.2 of the main paper. The spread estimates used alongside
them are proved next.

\section{Uniform spread around the mean}\label{supp:boundary}
\subsection{Spread in the constant-service model}
Here time is rescaled so that service capacity is one, and $\rho\in[1/2,1)$. Write $M=N^{-1}\sum_iQ_i$.

\begin{proof}[Proof of Lemma~3.3 of the main paper: constant service]\label{prf:lem:boundary-spread}
For $t\in(0,1]$ put $y_i=Q_i-M$,
$\Phi=\sum_i(e^{ty_i}+e^{-ty_i})$ and $f_i=2\sinh(ty_i)$.
We need the following deterministic inequality. If
$y_1\le\cdots\le y_N$ and $\sum_i y_i=0$, define
\[
 w_i=\frac{N+1-2i}{N},\qquad
 D=-\sum_iw_if_i=\frac1N\sum_{i<j}(f_j-f_i).
\]
Then
\begin{equation}\label{eq:boundary-rank}
 \Phi\le2N+10D.
\end{equation}
Let $Y$ be uniform on the $y_i$, $m$ a median, and
$H=\E|2\sinh(tY)-2\sinh(tm)|$. Median optimality gives
$\E|2\sinh(tY)-2\sinh(tY')|\ge H$ for an independent $Y'$, hence
$D\ge NH/2$.
If $m>0$, then $\E Y^- =\E Y^+\ge m/2$ and
$p=\PP(Y<0)\in(0,1/2]$. The function $f(y)=2\sinh(ty)$ is convex
on $[0,\infty)$ and $f(0)=0$. Conditional Jensen's inequality and convexity give
\[
 H\ge\E[f(-Y)\one_{\{Y<0\}}]
 \ge p f\!\left(\frac{\E Y^-}{p}\right)
 \ge p f\!\left(\frac{m}{2p}\right)
 \ge\tfrac12 f(m).
\]
For the last inequality, use $f(rx)\ge r f(x)$ for $r\ge1$.
The case $m<0$ follows by replacing $Y$ by $-Y$, and $m=0$ is
immediate. Thus $|f(m)|\le2H$ and $\E|f(Y)|\le3H$.
Together with $D\ge NH/2$ and $2\cosh x\le2+2|\sinh x|$,
these bounds prove~\eqref{eq:boundary-rank}.

Let $b_i$ be the arrival rates and $s_i=\one_{\{Q_i>0\}}$ the
service rates. An arrival at $j$ changes $y_i$ by
$\one_{\{i=j\}}-N^{-1}$, and a departure changes it by the negative
of this vector. The first-order arrival drift of $\Phi$ is
\[
 t\sum_i(b_i-\rho)f_i\le t\rho\sum_iw_if_i=-t\rho D.
\]
The minimum of $d\ge2$ samples is stochastically at most the minimum
of two, whose rank-$i$ probability is $[2(N-i)+1]/N^2$.
Monotonicity of $f_i$ gives the inequality; tied ranks have equal $f_i$.
The first-order service drift is
\[
 -t\sum_i\left(s_i-N^{-1}\sum_js_j\right)f_i\le0,
\]
because $s_i$ and $f_i$ are increasing functions of $Q_i$.

For $|h|\le1$, the exponential remainder is at most
$\tfrac12t^2e^t h^2$ times the corresponding exponential.
For an event at server $j$, the sum of the squared-coordinate
remainders is therefore bounded by
\(C t^2\left(e^{ty_j}+e^{-ty_j}+\Phi/N^2\right).\)
Using $b_j\le d\rho$, $s_j\le1$, and total queue-jump rate at most
$2N$, all second-order terms are bounded by $C'_dt^2\Phi$.
By~\eqref{eq:boundary-rank},
\(\mathcal L\Phi \le -t\rho\Phi/10+t\rho N/5+C'_dt^2\Phi.\)
Choose $t=t_d>0$ small enough, using $\rho\ge1/2$, to obtain
$\mathcal L\Phi\le-c_d\Phi+C''_dN$ with $c_d>0$.

To establish integrability, use the bounded concave test
$\psi_K(v)=K(1-e^{-v/K})$. Since
$\mathcal L\psi_K(\Phi)\le e^{-\Phi/K}\mathcal L\Phi$
and the jump rate is bounded for fixed $N$, stationarity yields
\(c_d\E[\Phi e^{-\Phi/K}] \le C''_dN\E e^{-\Phi/K}\le C''_dN.\)
Monotone convergence as $K\to\infty$ proves $\E\Phi\le C_dN$,
which implies~(27) of the main paper.
\end{proof}

\subsection{Spread with independent service environments}

\begin{proof}[Proof of Lemma~3.3 of the main paper: environments]\label{prf:lem:critical-mod-spread}
Set $y_i=Q_i-M_N$, $\phi_i=e^{ty_i}+e^{-ty_i}$ and $\Phi=\sum_i\phi_i$.
The constant-service argument above, with time multiplied by $\bc$,
gives $\mathcal L^0\Phi\le-c_0\Phi+C_0N$ for a fixed small
$t>0$. Here $\mathcal L^0$ replaces every capacity by $\bc$.
For a queue-only function $F$ define
\[
 S_i^F(x)=\one_{\{x_i>0\}}[F(x-e_i)-F(x)],
 \qquad B_F(x,e)=\sum_i h_{e_i}S_i^F(x).
\]
The Poisson equation~(21) of the main paper gives the exact identity
\begin{equation}\label{eq:supp-spread-corrected-test}
 \mathcal L(F+B_F/\kappa)
       =\mathcal L^0F+\kappa^{-1}\mathcal L_{\rm q}B_F.
\end{equation}
For $F=\Phi$, we claim
\begin{equation}\label{eq:critical-spread-correction}
 |B_\Phi|\le C\Phi,\qquad
 |\mathcal L_{\rm q}B_\Phi|\le C\Phi .
\end{equation}
Put $a_{il}=\one_{\{i=l\}}-N^{-1}$. A first difference of $\Phi$
at $i$ is bounded by $C(\phi_i+\Phi/N)$, proving the first bound.
For a jump at $j$, the $i=j$ summand of $B_\Phi$, including its busy
indicator, changes by at most $C(\phi_j+\Phi/N)$. For $i\ne j$ the
indicator is fixed and the exponential mixed difference is at most
$C\sum_l\phi_l|a_{il}a_{jl}|$. Its sum over $i\ne j$ is
$C(\phi_j+\Phi/N)$: coordinates $l=i,j$ contribute
$C\Phi/N+C\phi_j$, and the rest $C\Phi/N$. Intermediate exponentials
change by a bounded factor because $|a_{il}|\le1$. Each server's
queue-jump rate is at most $d\bc+c_*$, so summation over $j$ proves
\eqref{eq:critical-spread-correction}.

For sufficiently large fixed $\kappa_0$,
$F_\kappa=\Phi+B_\Phi/\kappa$ lies between $\Phi/2$ and $3\Phi/2$,
and~\eqref{eq:supp-spread-corrected-test} gives
$\mathcal LF_\kappa\le-c\Phi+CN$ with $c>0$ independent of $N,\rho,\kappa$.
Apply stationarity to the bounded increasing concave function
$\psi_K(F_\kappa)=K(1-e^{-F_\kappa/K})$. Its tangent inequality gives
\(c\E[\Phi e^{-F_\kappa/K}]\le CN.\)
The bounded jump rate justifies these tests. As $K\to\infty$,
monotone convergence gives the bound without any prior moment assumption.
\end{proof}

\section{Coupling and conditioned stationary references}\label{supp:coupling}
\subsection{Wasserstein preparation}

\begin{proof}[Proof of Lemma~4.5 of the main paper]\label{prf:lem:nonstationary-W1-input}
For $x\le y$, monotonicity of $\omega$ orders $\Pi_x$ below $\Pi_y$.
The quantile coupling attains their mean difference $y-x$, a lower
bound for every coupling, proving the isometry.

Since $\omega(x+1)=\omega(x)/d$, $\Pi_{x+1}$ translates $\Pi_x$ by one.
Reduce $x$ modulo one; the parameters $\omega(x)$ then form a compact
subset of $(0,\infty)$, giving uniform doubly exponential right tails and
\(\Pi_x((-\infty,-j]) =1-e^{-\omega(x)d^{-j+1}}\le C d^{-j},\qquad j\ge1,\)
uniformly for $x\in[0,1]$. Consequently, for any fixed sufficiently
small $\eta\in(0,\log d)$,
\begin{equation}\label{eq:nonstationary-reference-exp-moment}
 \sup_{x\in\mathbb R}
   \int e^{\eta|z-x|}\,\Pi_x(\dd z)<\infty.
\end{equation}
Reduce $\eta$ further if necessary so that
the applicable assertion of Lemma~3.3 of the main paper also gives
\[
 \E\int e^{\eta|z-X_N|}\,\mu_N^q(\dd z)
 =\E\frac1N\sum_i e^{\eta|Q_i-N^{-1}\sum_lQ_l|}\le C.
\]

For any measures $\mu,\nu$ with finite first moments, any real center
$x$, and $R>0$, clipping to $[x-R,x+R]$ and using its diameter yields
\begin{align*}
 W_1(\mu,\nu)
 &\le 2R\,d_{\TV}(\mu,\nu)\\
 &\quad+\int |z-x|\one_{\{|z-x|>R\}}\,\mu(\dd z)
       +\int |z-x|\one_{\{|z-x|>R\}}\,\nu(\dd z).
\end{align*}
Clipping into the real line does not increase total variation.
With $x=X_N$, the two exponential moment bounds give
\(\E W_1(\mu_N^q,\Pi_{X_N}) \le 2R\,\E d_{\TV}(\mu_N^q,\Pi_{X_N})+C e^{-\eta R/2}.\)
Use~(53) of the main paper, then let $N\to\infty$ and
$R\to\infty$. No moment bound on $X_N$ is required.
\end{proof}

\subsection{Arbitrary initial phases}

\begin{proof}[Proof of Lemma~4.6 of the main paper]\label{prf:lem:nonstationary-path-coupling}
Pairing ordered queue lengths realizes the initial Wasserstein
distance; permute phases with their queues. Use common arrivals and
sampled lists. Give distinct sampled labels common fresh independent
continuous priorities, using a single mark for repeated occurrences.
Among the shortest sampled queues, selecting the one with the lowest mark
implements uniform tie-breaking among distinct minima, including
repeated-index lists.
By symmetry among queues at a fixed length, this convention gives
exactly the birth rates~(2) of the main paper.

For prearrival vectors $q,r$, let $i,j$ be the chosen indices.
If $i=j$, $\ell_1$ distance is unchanged. Otherwise an increase
would be by two, requiring $q_i\ge r_i$ and $r_j\ge q_j$. Minimality gives
\(q_i\le q_j\le r_j\le r_i\le q_i.\)
Thus both are minima in both copies, contradicting their common
strict priority order. Every arrival is therefore nonexpansive.

At each pair, run the phases independently until their meeting time
$\tau_i$, then use a common trajectory, independently across pairs.
Both marginals retain generator $\kappa A$, and
\begin{equation}\label{eq:nonstationary-phase-meeting}
 \sup_{e,f}\E_{e,f}\tau_i\le C_A/\kappa.
\end{equation}
Indeed, irreducibility gives positive entries in $P=e^{s_0A}$ for
$s_0>0$. Two independent unit-speed copies agree after $s_0$ with
probability at least $\min_{e,f}\sum_gP(e,g)P(f,g)>0$.
Iteration bounds the meeting time by a geometric number of such
blocks; rescaling gives~\eqref{eq:nonstationary-phase-meeting}.

At each pair use a common rate-$c_*$ potential-service clock and shared
independent uniform marks, independently across pairs and of the phase
coupling. A copy accepts a mark below its capacity divided by $c_*$
when its queue is nonempty. After $\tau_i$, services are nonexpansive.
Before $\tau_i$, each potential-service event can increase the distance
by at most one. Let $B_N$ count potential services before the respective
meeting times.
 Then
\[
 \E B_N=c_*\sum_i\E\tau_i\le CN/\kappa.
\]
This proves the all-time bound, also at zero capacity and empty queues.
For constant service, common clocks are nonexpansive from time zero,
so $B_N=0$.

Since $\Delta_N(s)\le N^{-1}\sum_i\one_{\{s<\tau_i\}}$,
\[
 \int_0^T\Delta_N(Nt/\bc)\,\dd t
 \le\frac{\bc}{N^2}\sum_i\tau_i.
\]
Taking expectations proves~(75) of the main paper.
Pairing the empirical atoms proves~(76) of the main paper:
an atom can disagree only if its queue length or its phase disagrees,
and a nonzero difference of integer queue lengths is at least one.

By permutation equivariance, every permutation-invariant observable
has the same path law from $z$ and $\sigma z$. This statewise identity
also permits a permutation measurable in both initial configurations.
Hence optimal relabeling preserves the conditioned reference mean and
empirical-measure path laws. It need not preserve the labeled initial
law; only these invariant observables are used below.
\end{proof}

\subsection{Conditioning a stationary reference}

\begin{proof}[Proof of Lemma~4.7 of the main paper]\label{prf:lem:nonstationary-conditioned-reference}
Since $p_\theta$ is positive and gives zero mass to the endpoints of
$I$, $\PP(X_N\in I)\to\int_Ip_\theta>0$. Conditioning weights the
stationary path law by
$\one_{\{\mathcal X_N(0)\in I\}}/\PP(X_N\in I)$.
Time-zero evaluation is $J_1$-continuous, and the indicator is
discontinuous only at a limiting null set. Thus
Proposition~4.1 of the main paper gives the conditioned path limit.
The limiting weight depends only on $X(0)$, preserving its independence
of the Brownian motion.

For any nonnegative path functional its conditioned expectation is
at most its stationary expectation divided by $\PP(X_N\in I)$.
Lemma~4.5 of the main paper gives
(77) of the main paper; stationarity gives
(78) of the main paper. The stationary profile assertion in
Proposition~4.1 of the main paper gives
(79) of the main paper by the same comparison.
\end{proof}

\section{Whole-line dynamics and finite-system approximation}\label{supp:initial-layer}
\subsection{Construction, contraction and conservation}

\begin{proof}[Proof of Lemma~5.1 of the main paper]\label{prf:lem:initial-layer-semiflow}
Let $v_j^*=\one_{\{j\le0\}}$ be the tails of a point mass at zero.
The tails of $\mathcal P_1(\mathbb Z)$ form a closed subset of
$v^*+\ell^1(\mathbb Z)$, and
\begin{equation}\label{eq:supp-tail-W1}
 W_1(\mu,\nu)=\sum_{j\in\mathbb Z}|v_j(\mu)-v_j(\nu)|,
 \qquad
 \int z\,\mu(\dd z)=\sum_j(v_j(\mu)-v_j^*).
\end{equation}
Quantile coupling and the first-moment tail formula give these identities.
Replace $z^d$ by $(0\vee(z\wedge1))^d$. The resulting field is globally
$2\bc(d+1)$-Lipschitz on $v^*+\ell^1$ and takes $v^*$ into $\ell^1$,
so successive approximation gives a unique global solution.
Its bounded difference system has nonnegative off-diagonal coefficients;
integrating-factor iteration preserves positivity in $\ell^\infty$.
Comparison with zero, one and a shifted copy preserves
$1\ge v_j\ge v_{j+1}\ge0$. The affine $\ell^1$ condition gives endpoint
limits one and zero, so the solution is a $\mathcal P_1$ tail and the
modification is inactive.

For two such solutions write $w_j=v_j-\widetilde v_j$. Exactly,
\begin{equation}\label{eq:initial-layer-difference-generator}
 \dot w_j=\bc\{a_{j-1}(s)w_{j-1}-(a_j(s)+1)w_j+w_{j+1}\},
 \qquad
 a_j(s)=\sum_{r=0}^{d-1}v_j(s)^{d-1-r}\widetilde v_j(s)^r\in[0,d].
\end{equation}
This is the forward equation of a birth--death chain with bounded
upward rates $\bc a_j(s)$ and downward rate $\bc$. Its evolution
$P(s,r)$ is positive by integrating-factor iteration and mass preserving
by absolute summation of the equations. Thus
$\|P(s,0)w\|_1\le\|w\|_1$, proving the contraction.
To conserve the mean, subtract the field at $v^*$: the difference has
the conservative form~\eqref{eq:initial-layer-difference-generator}
and sums absolutely to zero. The field at $v^*$ has two opposite
nonzero entries, so~\eqref{eq:supp-tail-W1} gives constant mean.

For strictness, if $w(0)\ne0$ has sum zero, choose $i_+,i_-$ with
positive and negative entries. Downward-only paths give, for
$k\le\min(i_+,i_-)$,
\begin{equation}\label{eq:initial-layer-kernel-overlap}
 P_{ki}(s,0)\ge
 e^{-\bc(d+1)s}\frac{(\bc s)^{i-k}}{(i-k)!}>0,
 \qquad i\in\{i_+,i_-\},\quad s>0.
\end{equation}
The rate product $\bc^{i-k}$, bounded exit rate and jump-time simplex
volume $s^{i-k}/(i-k)!$ give this bound. Hence
$p=P(s,0)w(0)^+$ and $n=P(s,0)w(0)^-$ overlap at $k$, and
\[
 \|w(s)\|_1
   =\|w(0)\|_1-2\sum_j\min(p_j,n_j)<\|w(0)\|_1.
\]
Equal initial means give the required zero sum of $w(0)$.

At equilibrium $v_j^d-v_{j+1}$ is constant and the right-tail limit
makes it zero. The endpoint limits force $0<v_j<1$, hence
$v_j=e^{-\omega d^j}$ for a unique $\omega>0$. The strictly decreasing bijection
$r$ then identifies $\omega(x)$ at prescribed mean $x$.
\end{proof}

\subsection{Approximation during the initial layer}
Use physical time throughout this subsection. The following quantitative
estimate supplies the compact-family approximation in
Lemma~5.3 of the main paper.

\begin{proposition}[Finite-system approximation bound]\label{prop:supp-fluid-error}
Let $N\ge2$, $\rho_N=1-N^{-1}$ and $m_N=\lfloor\log_dN\rfloor$.
For a deterministic initial queue--phase configuration and $\eta>0$, put
\[
 K_N=N^{-1}\sum_i e^{\eta|Q_i(0)-m_N|},\qquad
 v_\eta=(d\bc+c_*)(e^\eta-1),\qquad L_F=2\bc(d+1).
\]
For every $T>0$ and $\kappa\ge1$, the whole-line semiflow $S_s$ satisfies
\begin{align}
 \E\sup_{0\le s\le T}W_1(P_N(s),S_sP_N(0))
 &\le C_{\eta}e^{L_FT}\bigg\{
       \frac{T}{N}+\frac{1+T}{\kappa}\notag\\
 &\qquad+\sqrt{\frac{K_NT}{N}}e^{v_{\eta}T/2}
              \left(1+\frac{\sqrt{\kappa+1}}{\kappa}\right)\notag\\
 &\qquad+K_NT e^{v_{\eta}T-\eta m_N}\bigg\}.
                                      \label{eq:initial-layer-fluid-error}
\end{align}
The constant is independent of $N,T$ and of the initial phases.
For constant service, omit the terms containing $\kappa$.
\end{proposition}

\begin{proof}\label{prf:prop:supp-fluid-error}
Set $Z_i(s)=Q_i(s)-m_N$, $\ell_N=-m_N$, and, for every integer $j$,
\[
 V_j^N(s)=\frac1N\sum_i\one_{\{Z_i(s)\ge j\}},\qquad
 p_j(s)=V_j^N(s)-V_{j+1}^N(s).
\]
Thus $V_j^N=1$ for $j\le\ell_N$. For a probability tail $v$ on the
whole line define
\(F_j(v)=\bc(v_{j-1}^d-v_j^d-v_j+v_{j+1}).\)
For two such tails with finite first moments,
\begin{equation}\label{eq:initial-layer-F-Lipschitz}
 \sum_j|F_j(v)-F_j(\widetilde v)|\le L_F\sum_j|v_j-\widetilde v_j|,
 \qquad W_1(\mu,\nu)=\sum_j|v_j-\widetilde v_j|.
\end{equation}
This follows from $|v^d-w^d|\le d|v-w|$ on $[0,1]$.

Each coordinate has queue-jump intensity at most $d\bc+c_*$. Its
generator applied to $e^{\eta|Z_i|}$ gives
\begin{equation}\label{eq:initial-layer-moment-growth}
 \E\sum_j e^{\eta|j|}p_j(s)
       \le K_N e^{v_{\eta}s},\qquad s\ge0.
\end{equation}
This uses no initial phase assumption. For $g_e=c_e-\bc$ and the
bounded solution $-Ah=\boldsymbol g$ from~(21) of the main paper, define
\[
 H_j(s)=\one_{\{j>\ell_N\}}\frac1N
             \sum_i h_{J_i(s)}\one_{\{Z_i(s)=j\}},\qquad
 G_j(s)=\one_{\{j>\ell_N\}}\frac1N
             \sum_i g_{J_i(s)}\one_{\{Z_i(s)=j\}}.
\]
Writing $L_q$ for the queue generator, the phase drift of $H_j$ is
$-\kappa G_j$, so Dynkin's formula gives
\begin{equation}\label{eq:initial-layer-phase-corrector}
 \int_0^sG_j(r)\,\dd r
 =\frac{H_j(0)-H_j(s)}{\kappa}
   +\frac1\kappa\int_0^s L_qH_j(r)\,\dd r
   +\frac{M_j^H(s)}{\kappa}.
\end{equation}
Pathwise $\sum_j|H_j|\le\|h\|_\infty$. Each queue jump changes
at most two masses, including at the empty boundary, so
\begin{equation}\label{eq:initial-layer-H-drift}
 \sum_j|L_qH_j|\le 2\|h\|_\infty(d\bc+c_*).
\end{equation}
The empirical tail equation is
\begin{align}
 V_j^N(s)&=V_j^N(0)+\int_0^sF_j(V^N(r))\,\dd r+M_j^q(s)\notag\\
   &\quad-\frac{\bc}{N}\int_0^s
                    (V_{j-1}^N(r)^d-V_j^N(r)^d)\,\dd r
     -\int_0^sG_j(r)\,\dd r \notag\\
   &\quad+\one_{\{j=\ell_N\}}\bc\int_0^sp_{\ell_N}(r)\,\dd r.
                                      \label{eq:initial-layer-tail-equation}
\end{align}
Above $\ell_N$ this is ordinary level crossing. At $\ell_N$, the
last term cancels the fictitious whole-line flux $-\bc p_{\ell_N}$;
other terms vanish, as do all terms below $\ell_N$. The nonnegative
differences $v_{j-1}^d-v_j^d$ sum to one, so the arrival correction
has $\ell_1$ norm at most $\bc T/N$.
Equation~\eqref{eq:initial-layer-moment-growth} bounds the expected
integrated boundary correction by
\begin{equation}\label{eq:initial-layer-boundary-error}
 \bc K_NT e^{v_{\eta}T-\eta m_N}.
\end{equation}

The predictable brackets satisfy
\begin{align*}
 \E\langle M_j^q\rangle_T
   &\le\frac C N\int_0^T\E(p_{j-1}+p_j)(s)\,\dd s,\\
 \E\langle M_j^H\rangle_T
   &\le\frac{C(\kappa+1)}N
              \int_0^T\E(p_{j-1}+p_j+p_{j+1})(s)\,\dd s.
\end{align*}
Upward and downward crossings start at $j-1$ and $j$, respectively,
with bounded per-queue rates. Phase jumps contribute at most $C\kappa p_j/N$;
queue jumps affect adjacent masses. Absent lower-boundary terms only
improve the bounds; queue--phase independence is not used.

If $A_j=\int_0^T\E(p_{j-1}+p_j+p_{j+1})(s)\,\dd s$, then
\[
 \sum_j\sqrt{A_j}
 \le\left(\sum_j e^{-\eta|j|}\right)^{1/2}
       \left(\sum_j e^{\eta|j|}A_j\right)^{1/2}
 \le C_\eta\sqrt{K_NT}\,e^{v_{\eta}T/2}.
\]
Doob's inequality, followed by Tonelli's theorem, therefore yields
\begin{align}
 \E\sum_j\sup_{s\le T}|M_j^q(s)|
       &\le C_\eta\sqrt{K_NT/N}\,e^{v_{\eta}T/2},\notag\\
 \E\sum_j\sup_{s\le T}|M_j^H(s)|
       &\le C_\eta\sqrt{K_NT(\kappa+1)/N}\,e^{v_{\eta}T/2}.
                                      \label{eq:initial-layer-martingale-sum}
\end{align}
Localization to finite index sets and these finite bounds justify
summing the tail equations. Subtract the equation for $S_sP_N(0)$,
substitute~\eqref{eq:initial-layer-phase-corrector}, and take the
supremum of the $\ell_1$ difference. The endpoint correction uses
$\sup_{s\le T}\sum_j|H_j(s)-H_j(0)|\le2\|h\|_\infty$, rather than
coordinatewise endpoint suprema. Gronwall with
\eqref{eq:initial-layer-F-Lipschitz} gives
\eqref{eq:initial-layer-fluid-error}. Constant $c$ gives $\boldsymbol g=h=0$.

\end{proof}

\begin{proof}[Proof of Lemma~5.3 of the main paper]\label{prf:lem:initial-layer-fluid}
For compact initial families, clip centered queues to $[-R,R]$, with
integer $R>0$, keeping phases. The uncentered queues are nonnegative
for large $N$, and this monotone pairing realizes
\[
 e_R(P_N(0)):=\int (|z|-R)_+\,P_N(0)(\dd z).
\]
Common phases, arrival choices and service marks make the coupling of
Lemma~4.6 of the main paper nonexpansive without a
discrepancy term. Whole-line contraction gives the same bound for the
fluid solutions, so clipping costs at most $2e_R(P_N(0))$.
For clipped data $K_N\le e^{\eta R}$, and
\eqref{eq:initial-layer-fluid-error} vanishes at $T=b_N=o(\log N)$:
$e^{Cb_N}=N^{o(1)}$, $m_N=\log_dN+O(1)$ and eventually $\kappa_N\ge N$.
The original error therefore has limsup at most
$2\sup_{P\in\mathcal K}e_R(P)\to0$ as $R\to\infty$, proving
(95) of the main paper.

For random starts, condition on the full initial configuration and
localize $P_N(0)$ to a compact set of arbitrarily high probability.
The uniform bound and Markov's inequality give convergence in probability
without moments on the exceptional set. Mean conservation and the
$W_1$ bound on mean differences give~(96) of the main paper.
\end{proof}

\section{An equivalent logistic coordinate}\label{sec:logistic-coordinate}
Write $\Pi_x([j,\infty))=e^{-\omega(x)d^j}$ for the mean-$x$ family
of the main paper. The limiting mean satisfies
$\dd X=(g_\theta(X)-1)\,\dd t+\sqrt2\,\dd B$,
where $g_\theta(x)=(d-1)d^\theta\omega(x)$ and
$\omega(x+1)=\omega(x)/d$.
The following logistic-type representation is exact. Put $\ell_d=\log d$
and $\chi_d(x)=d^x\omega(x)$, so that $\chi_d(x+1)=\chi_d(x)$ by~(6) of the main paper.
For $R_t=d^{-X(t)}>0$, It\^o's formula, with $W=-B$, gives
\[
 \dd R_t=\ell_d R_t\bigl[1+\ell_d-(d-1)d^\theta R_t\chi_d(-\log_d R_t)\bigr]\,\dd t
                  +\sqrt2\,\ell_d R_t\,\dd W_t.
\]
The quadratic coefficient is periodic in $\log_d R_t$.
The correction $\ell_d^2R_t\,\dd t$ comes from $\dd[X]_t=2\,\dd t$.
For the accelerated centered empirical mean $\mathcal X_N$ defined in
the main paper, continuous mapping and Theorem~2.2 of the main paper give
$d^{-\mathcal X_N}\Rightarrow R$ in the $J_1$ topology on finite time intervals.


\begin{thebibliography}{99}
\bibitem[Atar and Wolansky(2025)]{AW}Atar, R., Wolansky, G.:
Invariance principle and McKean--Vlasov limit for randomized load balancing in heavy traffic.
Ann. Appl. Probab. \textbf{35}(5), 3046--3085 (2025).
\href{https://doi.org/10.1214/25-AAP2145}{doi:10.1214/25-AAP2145}. \MR{4975042}
\bibitem[Banerjee et~al.(2026)Banerjee, Budhiraja, and Loeser]{BBL2026}
Banerjee, S., Budhiraja, A., Loeser, E.:
Rank based routing in large server systems under extreme congestion.
Preprint (2026). \ARXIV{2605.17677}.
\bibitem[Bhamidi et~al.(2022)Bhamidi, Budhiraja, and Dewaskar]{BBD2022}
Bhamidi, S., Budhiraja, A., Dewaskar, M.:
Near equilibrium fluctuations for supermarket models with growing choices.
Ann. Appl. Probab. \textbf{32}(3), 2083--2138 (2022).
\href{https://doi.org/10.1214/21-AAP1729}{doi:10.1214/21-AAP1729}.
\bibitem[Bordenave et~al.(2010)Bordenave, McDonald, and Proutiere]{BMP}Bordenave, C., McDonald, D.R., Proutiere, A.:
A particle system in interaction with a rapidly varying environment: mean field limits and applications.
Netw. Heterog. Media \textbf{5}(1), 31--62 (2010).
\href{https://doi.org/10.3934/nhm.2010.5.31}{doi:10.3934/nhm.2010.5.31}. \MR{2601988}
\bibitem[Brightwell et~al.(2018)Brightwell, Fairthorne, and Luczak]{BFL2018}
Brightwell, G., Fairthorne, M., Luczak, M.J.:
The supermarket model with bounded queue lengths in equilibrium.
J. Stat. Phys. \textbf{173}, 1149--1194 (2018).
\href{https://doi.org/10.1007/s10955-018-2044-7}{doi:10.1007/s10955-018-2044-7}.
\bibitem[Budhiraja and Friedlander(2019)]{BF2019}
Budhiraja, A., Friedlander, E.:
Diffusion approximations for load balancing mechanisms in cloud storage systems.
Adv. Appl. Probab. \textbf{51}(1), 41--86 (2019).
\href{https://doi.org/10.1017/apr.2019.3}{doi:10.1017/apr.2019.3}.
\bibitem[Eschenfeldt and Gamarnik(2016)]{EG2016}Eschenfeldt, P., Gamarnik, D.:
Supermarket queueing system in the heavy traffic regime. Short queue dynamics.
Preprint (2016). \ARXIV{1610.03522}.
\bibitem[Gast(2017)]{Gast2017}Gast, N.:
Expected values estimated via mean-field approximation are $1/N$-accurate.
Proc. ACM Meas. Anal. Comput. Syst.
\textbf{1}(1), 1--26 (2017).
\href{https://doi.org/10.1145/3084454}{doi:10.1145/3084454}.
\bibitem[Graham(2005)]{Graham2005}
Graham, C.:
Functional central limit theorems for a large network in which customers
join the shortest of several queues.
Probab. Theory Related Fields \textbf{131}, 97--120 (2005).
\href{https://doi.org/10.1007/s00440-004-0372-9}{doi:10.1007/s00440-004-0372-9}.
\bibitem[Gupta and Walton(2019)]{GW2019}
Gupta, V., Walton, N.:
Load balancing in the nondegenerate slowdown regime.
Oper. Res. \textbf{67}(1), 281--294 (2019).
\href{https://doi.org/10.1287/opre.2018.1768}{doi:10.1287/opre.2018.1768}.
\bibitem[Hurtado-Lange and Grosof(2026)]{HLG2026}Hurtado-Lange, D., Grosof, I.:
Markov modulated JSQ in heavy traffic via the Poisson equation.
Preprint (2026). \ARXIV{2603.00387}.
\bibitem[Hurtado-Lange and Maguluri(2022)]{HLM2022}
Hurtado-Lange, D., Maguluri, S.T.:
A load balancing system in the many-server heavy-traffic asymptotics.
Queueing Syst. \textbf{101}, 353--391 (2022).
\href{https://doi.org/10.1007/s11134-022-09847-7}{doi:10.1007/s11134-022-09847-7}.
\bibitem[Jaleel et~al.(2026)Jaleel, Delasay, and Doroudi]{JDD2026}Jaleel, J.A., Delasay, M., Doroudi, S.:
Scalable load balancing in interference-prone queueing systems.
Ann. Oper. Res. \textbf{361}(2), 637--679 (2026).
\href{https://doi.org/10.1007/s10479-026-07217-7}{doi:10.1007/s10479-026-07217-7}.
\bibitem[Janson(2006)]{Janson2006}Janson, S.:
Rounding of continuous random variables and oscillatory asymptotics.
Ann. Probab. \textbf{34}(5), 1807--1826 (2006).
\href{https://doi.org/10.1214/009117906000000232}{doi:10.1214/009117906000000232}. \MR{2271483}
\bibitem[Luczak and McDiarmid(2006)]{LM2006}
Luczak, M.J., McDiarmid, C.:
On the maximum queue length in the supermarket model.
Ann. Probab. \textbf{34}(2), 493--527 (2006).
\href{https://doi.org/10.1214/00911790500000710}{doi:10.1214/00911790500000710}.
\bibitem[Luczak and McDiarmid(2007)]{LM2007}
Luczak, M.J., McDiarmid, C.:
Asymptotic distributions and chaos for the supermarket model.
Electron. J. Probab. \textbf{12}, 75--99 (2007).
\href{https://doi.org/10.1214/EJP.v12-391}{doi:10.1214/EJP.v12-391}.
\bibitem[Luczak and Norris(2013)]{LN}Luczak, M.J., Norris, J.R.:
Averaging over fast variables in the fluid limit for Markov chains: application to the supermarket model with memory.
Ann. Appl. Probab. \textbf{23}(3), 957--986 (2013).
\href{https://doi.org/10.1214/12-AAP861}{doi:10.1214/12-AAP861}. \MR{3076675}
\bibitem[Mitzenmacher(2001)]{Mitzenmacher2001}Mitzenmacher, M.:
The power of two choices in randomized load balancing.
IEEE Trans. Parallel Distrib. Syst. \textbf{12}(10), 1094--1104 (2001).
\href{https://doi.org/10.1109/71.963420}{doi:10.1109/71.963420}.
\bibitem[Oseledets and Khmelev(2000)]{OK2000}
Oseledets, V.I., Khmelev, D.V.:
Global stability of infinite systems of nonlinear differential equations
and nonhomogeneous countable Markov chains.
Problems Inform. Transmission \textbf{36}(1), 54--70 (2000).
Russian original: Probl. Peredachi Inf. \textbf{36}(1), 60--76.
\url{https://www.mathnet.ru/eng/ppi470}.
\bibitem[Varma et~al.(2025)Varma, Castro, and Maguluri]{VCM2025}
Varma, S.M., Castro, F., Maguluri, S.T.:
Power-of-$d$ choices load balancing in the sub-Halfin--Whitt regime.
Stoch. Syst. \textbf{15}(4), 291--336 (2025).
\href{https://doi.org/10.1287/stsy.2024.0063}{doi:10.1287/stsy.2024.0063}.
\bibitem[Vvedenskaya and Suhov(2005)]{VS2005}Vvedenskaya, N.D., Suhov, Y.M.:
Dynamic routing queueing systems with vacations.
Inform. Protsessy \textbf{5}(1), 74--86 (2005).
\bibitem[Vvedenskaya et~al.(1996)Vvedenskaya, Dobrushin, and Karpelevich]{VDK}Vvedenskaya, N.D., Dobrushin, R.L., Karpelevich, F.I.:
Queueing system with selection of the shortest of two queues: an asymptotic approach.
Probl. Peredachi Inf. \textbf{32}(1), 20--34 (1996);
English translation: Problems Inform. Transmission \textbf{32}(1), 15--27.
\url{https://www.mathnet.ru/ppi298}. \MR{1384927}
\bibitem[Whitt(2007)]{Whitt2007}Whitt, W.:
Proofs of the martingale FCLT.
Probab. Surv. \textbf{4}, 268--302 (2007).
\href{https://doi.org/10.1214/07-PS122}{doi:10.1214/07-PS122}.
\bibitem[Ying(2017)]{Ying2017}Ying, L.:
Stein's method for mean-field approximations in light and heavy traffic regimes.
Proc. ACM Meas. Anal. Comput. Syst.
\textbf{1}(1), 1--27 (2017).
\href{https://doi.org/10.1145/3084449}{doi:10.1145/3084449}.
\bibitem[Kalimulina(2026)]{ASupp}
Kalimulina, E. Yu. (2026).
Supplement to ``Critical heavy traffic and diffusion limits for shortest-of-$d$ dynamic routing.''
\end{thebibliography}
\end{document}